\documentclass[11pt]{article}

\usepackage{amsmath, amsfonts, amssymb}
\usepackage{graphicx}
\usepackage{algorithmic}
\usepackage{algorithm}
\usepackage{bm}
\usepackage{xcolor}
\usepackage{booktabs}
\usepackage{siunitx}
\usepackage{caption}
\usepackage{makecell}
\usepackage{subcaption}
\usepackage{hyperref}
\usepackage{bibentry}
\usepackage{amsopn}
\usepackage{lineno}
\usepackage{float}
\usepackage{makecell}
\DeclareGraphicsExtensions{.eps,.pdf,.png,.jpg}

\usepackage{amsthm}
\newtheorem{theorem}{Theorem}[section]
\newtheorem{lemma}{Lemma}[section]
\newtheorem{corollary}{Corollary}[section]
\newtheorem{assumption}{Assumption}[section]
\newtheorem{remark}{Remark}[section]

\newtheorem{definition}{Definition}[section]
\makeatletter
\@ifundefined{proof}{%
  \newenvironment{proof}[1][Proof]{%
    \par\noindent\textit{#1.}\enspace
  }{%
    \hfill\ensuremath{\square}\par
  }%
}{}
\makeatother
\numberwithin{equation}{section}
\definecolor{dkgreen}{rgb}{0,0.6,0}

\newcommand{\R}{\mathbb{R}}
\newcommand{\E}{\mathbb{E}}

\newcommand{\be}{\begin{equation}}
\newcommand{\ee}{\end{equation}}
\newcommand{\bee}{\begin{equation*}}
\newcommand{\eee}{\end{equation*}}
\newcommand{\bea}{\begin{eqnarray}}
\newcommand{\eea}{\end{eqnarray}}
\newcommand{\beaa}{\begin{eqnarray*}}
\newcommand{\eeaa}{\end{eqnarray*}}
\allowdisplaybreaks[4]

\title{Mirror Descent Linearized Augmented Lagrangian Methods for Nonconvex Constrained Stochastic Zeroth-Order Optimization}
\author{
  Qiankun Shi\thanks{Sun Yat-sen University, Guangzhou, 510006, China; Pengcheng Laboratory, Shenzhen, 518066, China. (\texttt{shiqk@pcl.ac.cn}).}
   \and 
   Han Yuan\thanks{Sun Yat-sen University, Guangzhou, 518066, China.(\texttt{hyuan\_0604@163.com})}\and
   Xiao Wang\thanks{Sun Yat-sen University, Guangzhou, 510006, China. (\texttt{wangx936@mail.sysu.edu.cn}).} \and
   Hao Wang\thanks{ShanghaiTech University, Shanghai, 201210, China. (\texttt{wanghao1@shanghaitech.edu.cn})}
}

\hypersetup{
  pdftitle={Mirror Descent Linearized Augmented Lagrangian Methods for Nonconvex Constrained Stochastic Zeroth-Order  Optimization},
  pdfauthor={Qiankun Shi, Hao Wang, Xiao Wang, Han Yuan}
}

\begin{document}

\maketitle

% Abstract
\begin{abstract}
In this paper, we study nonconvex constrained stochastic zeroth-order optimization problems in which exact constraint information is available, while the objective can be accessed only through stochastic function evaluations. To solve this class of problems, we propose a framework of mirror descent linearized augmented Lagrangian methods that employs two-point stochastic zeroth-order gradient estimators and exploits non-Euclidean mirror descent geometry. Under stochastic oracle assumptions and a regularity assumption, we establish oracle complexity guarantees for finding an $\epsilon$-KKT point under mirror descent geometries parameterized by $p \geq 2$. Under Rademacher smoothing, our analysis reveals a trade-off between the dimension-dependent variance of the zeroth-order gradient estimators and the smoothness of the non-Euclidean mirror map. In the high accuracy regime, the resulting effective oracle complexity is $\mathcal{O}(p d^{2/p}\epsilon^{-3})$ for $p \in [2,2\ln d]$ and $\mathcal{O}(\ln d\,\epsilon^{-3})$ for $p > 2\ln d$. These bounds reduce the dependence on the problem dimension $d$ in the leading oracle complexity term at the expense of additional geometry dependent lower-order terms. In particular, when $p=2$, the method recovers the standard Euclidean setting and achieves an oracle complexity of $\mathcal{O}(d\epsilon^{-3})$, improving the dependence on $\epsilon$ over existing constrained stochastic zeroth-order methods. Furthermore, to eliminate near-feasibility requirement on the initial point, we develop a multi-stage scheme that finds an $\epsilon$-KKT point within $\mathcal{O} (1+\log\log(e/\epsilon))$ stages while preserving the leading-order high-accuracy oracle complexity, up to an additional $\epsilon$ independent initialization cost. Numerical tests on quadratically constrained quadratic programs, black-box adversarial attacks, and fairness-constrained classification showcase the effectiveness of the proposed method.
\end{abstract}

% Keywords
\textbf{Keywords}: Nonconvex constrained optimization, stochastic optimization, zeroth-order optimization, mirror descent, oracle complexity, high-accuracy regime

% MSC codes
\textbf{MSC codes}: 65K05, 68Q25, 90C30

\section{Introduction}
In this paper, we focus on the following constrained optimization problem:
\be\label{p}
\begin{aligned}
    \min_{x \in X}\quad  &\Phi(x):=\{f(x):= \E_\xi[F(x;\xi)]\} +  h(x)\quad \text{s.t.} \quad c(x) = 0,
\end{aligned}
\ee
where $X \subset \R^d$ is convex, $f: \R^d \rightarrow \R$ and $c: \R^d \rightarrow \R^m$ are continuously differentiable  and potentially  nonconvex, $\xi$ is a random variable from the probability space $(\Omega, \mathcal F, \mathbb P)$ and is independent of $x$, {$F(x;\xi)$} is continuously differentiable in  $x$ for almost any $\xi$, $\E_\xi[\cdot]$ represents the expectation taken with respect to (w.r.t.) $\xi$, and {$h:\R^d \to \R\cup \{\infty\} $} is a simple convex {but} nonsmooth function. 
Problem \eqref{p} arises in many application fields. For instance, to enforce specific behaviors or properties in deep learning, constraints are imposed on the output of deep neural networks \cite{conlab}, such as the physics-constrained deep learning model \cite{ZHU201956}, constraint-aware deep neural network compression \cite{Chen_2018_ECCV}, manifold-regularized deep learning \cite{roy2018}. 
Other applications include portfolio allocation \cite{bertsimas2018robust}, {two-stage and} multi-stage modeling \cite{shapiro2021lectures} and constrained maximum likelihood estimation \cite{chatterjee2016constrained}. 
Throughout this paper, we assume that $d\ge3$ and that for the objective function $f$ only stochastic function values are available at a query point.

Zeroth-order algorithms, also known as gradient-free algorithms, address optimization problems without relying on explicit gradient computations. These algorithms can be roughly  categorized into two classes: direct search methods \cite{hooke1961direct,nelder1965simplex,rosenbrock1960automatic} and model-based methods \cite{billups2013derivative,conn2008geometry2,powell2001lagrange}. The former class operates by sampling the objective function and selecting iterates through straightforward comparisons of function values, while the latter approximates the original optimization problem with simpler models and determines iterates based on these approximations. For broader surveys, the reader is referred to \cite{conn2009introduction,larson2019derivative}.  Among model-based methods, finite-difference methods have been studied comprehensively.  
For instance, Nesterov's seminal work \cite{Nesterov17} establishes the complexity bound for gradient-free methods in convex optimization. Duchi et al. \cite{duchi2015optimal} introduce two-point gradient estimators for stochastic convex minimization, achieving a nearly optimal oracle complexity of $\mathcal{O}(d\epsilon^{-2})$.  
In stochastic zeroth-order nonconvex optimization, general methods such as zeroth-order SGD \cite{sgd} exhibit an oracle complexity order of $\mathcal{O}(d\epsilon^{-4})$ and an improved order of   $\mathcal{O}(d\epsilon^{-3})$ when  incorporating variance reduction techniques \cite{spider, huang2020accelerated, JWZL}. 
A key challenge in this field is the impact of the problem dimension \(d\), especially in high-dimensional settings, as zeroth-order gradient estimators yield a
dimension-dependent variance, resulting in a strong dependence of oracle complexity on $d$ \cite{liu2020primer,ma2025revisiting}. To address this, stronger assumptions such as sparsity are often introduced \cite{BG22,roy2022stochastic,wang2018stochastic}. 
However, such sparsity assumptions may not always hold. 
Moreover, the aforementioned zeroth-order algorithms are designed for unconstrained settings, making them inapplicable to general constrained optimization problems. 
%To alleviate the polynomial dependence on the problem dimension $d$ in stochastic optimization, utilizing non-Euclidean geometries via mirror descent is an effective approach. 

For  constrained stochastic optimization, stochastic first- and second-order methods have been well developed.
Among stochastic first-order methods, 
proximal point methods including \cite{icpp,boob2022level,proxcm-ytb} focus on nonconvex inequality-constrained optimization and solve the original problem through a sequence of strongly convex subproblems. 
Stochastic sequential quadratic programming (SQP) methods have also recently been studied for equality constrained stochastic optimization~\cite{berahas2021sequential,COR,fang2024fully}, with a complexity analysis provided in \cite{COR}. In  subsequent work,  Curtis et al. \cite{Curtis23} study problems with both equality and inequality constraints and present a stochastic SQP method with global convergence analysis. 
Penalty methods for nonconvex constrained stochastic optimization discourage constraint violations by  incorporating a penalty term, weighted by a parameter \cite{alacaoglu2024complexity,cui2025exact,li2026mossp,spd,lu2026variance,SWW,wang2017penalty,pdsg}. 
In terms of stochastic second-order algorithms, Na et al. \cite{na2023adaptive} propose an adaptive SQP method with differentiable exact augmented Lagrangians, while Wang \cite{Wang2024Complexity} employs a stochastic cubic-regularized primal-dual method to obtain an approximate second-order stationary point. More recently, Shi et al.~\cite{shi2026curvature} proposed a
curvature-oriented variance-reduction method for finding approximate
second-order stationary points of unconstrained nonconvex stochastic
optimization. They further developed an augmented Lagrangian variant,
for problems with deterministic nonlinear equality
constraints. 
However, all these methods rely on unbiased stochastic gradients and are therefore inapplicable when only biased stochastic gradient estimates are available.

Although in the literature Frank-Wolfe methods and projected gradient methods based on zeroth-order estimation have been studied for set-constrained optimization \cite{BGNips18, BG22}, the nonconvexity of constraints in  nonconvex constrained optimization limits the study of zeroth-order algorithms. In \cite{wang2017penalty}, a double-loop penalty method, which calls a stochastic zeroth-order algorithm to solve each penalty subproblem, is proposed for equality constrained stochastic optimization, but with a relatively high oracle complexity.  
Recently, stochastic zeroth-order algorithms for solving functionally constrained optimization problems with stochastic inequality constraints have been studied in \cite{nguyen2022stochastic}, yielding an oracle complexity of order \(\mathcal{O}(d\epsilon^{-6})\).
Notably, both algorithms exhibit a complexity in $d$ of at least \(\mathcal{O}(d)\). 
Reformulating the original problem into a min-max optimization problem, Shi et al. present a doubly stochastic zeroth-order gradient algorithm in \cite{SGG} for heavily constrained nonconvex optimization with numerous constraints, achieving the complexity of \(\mathcal{O}(\epsilon^{-4})\). 
However, \cite{SGG} hides the dimensional dependence in its assumptions. In a related deterministic black-box setting, Jin et al.~\cite{jin2026query} consider nonconvex constrained problems where both the objective and constraint functions are accessed through zeroth-order oracles. Their block zeroth-order GDA-type algorithms reduce the number of queries needed for a single-step gradient estimate to \(\mathcal{O}(1)\), while maintaining the best-known overall query complexity \(\mathcal{O}(d\epsilon^{-4})\) for finding an \(\epsilon\)-stationary solution. More recently, Na~\cite{na2025derivative} developed a derivative-free
stochastic SQP method for problems with a stochastic objective and
deterministic equality constraints, both accessible only through
zeroth-order oracles, but without providing a complexity analysis. Alongside our investigation in the zeroth-order setting, recent work by Gonz\'{a}lez et al. \cite{gonzalez2026mirror} similarly applies stochastic mirror descent in the $\ell_1$ setup to achieve nearly dimension-independent convergence rates for differentially-private saddle-point problems. 
% \blue{Although non-Euclidean mirror geometries have been shown to improve 
% dimension dependence in stochastic first-order optimization by adapting 
% the algorithm to the underlying geometry of the problem 
% \cite{duchi2010composite,gonzalez2026mirror}, 
% their extension to constrained zeroth-order optimization remains nontrivial. The main difficulty lies in the fact that zeroth-order gradient estimators 
% are biased and suffer from dimension-dependent estimation errors, while the 
% augmented Lagrangian framework requires simultaneously controlling objective 
% descent and constraint feasibility through primal-dual updates.}

{\bf Contributions.} 
The main contributions of this work are summarized as follows.
\begin{itemize}
\item
We propose a single-loop mirror descent linearized augmented Lagrangian method for nonconvex constrained stochastic zeroth-order optimization. 
Since neither the exact objective value nor its gradient is available, we construct a two-point stochastic zeroth-order gradient estimator and incorporate a momentum-based variance-reduction mechanism. 
At each iteration, the primal variable is updated through a mirror descent proximal subproblem based on a stochastic linearization of the augmented Lagrangian, while the multiplier is updated by a single dual step. 
This framework allows the primal geometry to be chosen beyond the standard Euclidean setting.

\item
Under stochastic oracle assumptions and a regularity  constraint qualification, we establish stationarity and primal-feasibility bounds under general mirror descent geometries parameterized by $p\geq 2$, and derive oracle complexity guarantees for reaching an $\epsilon$-KKT point.
For $p\in[2,2\ln d]$, the stochastic estimation error is controlled directly in the $\ell_p$-norm.
For $p>2\ln d$, we retain the original mirror geometry associated with $p$, while controlling the stochastic estimation error in the smaller norm exponent
$\bar p=2\ln d$.
This separation prevents the stochastic-estimation constants from deteriorating as $p$ increases and provides a unified analysis for the full range $p\geq2$.

\item
% To realize the non-Euclidean mirror geometry required by our analysis, we construct the smoothed distance-generating function
% \[
% v_\delta(x)
% =
% \frac{1}{2(q-1)}
% (
% \sum_{i=1}^d (x_i^2+\delta^2)^{q/2}
% )^{2/q},
% \qquad
% \frac{1}{p}+\frac{1}{q}=1.
% \]
% The smoothing parameter $\delta>0$ removes the Hessian singularity of the standard squared $\ell_q$-norm near the origin. We prove that $v_\delta$ remains $1$-strong convexity w.r.t. the $\ell_q$-norm and that its gradient is $L_v$-Lipschitz continuous on the constraint set. 
Under Rademacher smoothing, we  make explicit the interaction between the dimension-dependent variance of the zeroth-order estimator and the smoothness of the proposed mirror map. For $p\in[2,2\ln d]$, the oracle complexity is
\(
\mathcal{O}(
p d^{2/p}\epsilon^{-3}
+p^2d\,\epsilon^{-2}
+p^3d^{2-2/p}
),
\) 
whereas for $p>2\ln d$ it is
\(
\mathcal{O}(
\ln d\cdot\epsilon^{-3}
+p d^{1-2/p}\ln d\cdot\epsilon^{-2}
+p^2d^{2-4/p}\ln d).
\) 
The first terms in these bounds reflect the dimension-dependent variance of zeroth-order gradient estimation, while the remaining terms arise from the smoothness of the mirror map. In the high-accuracy regime, the resulting effective oracle complexities reduce to
$\mathcal{O}(p d^{2/p}\epsilon^{-3})$ and
$\mathcal{O}(\ln d\cdot\epsilon^{-3})$, respectively.
In particular, when $p=2$, the proposed mirror geometry reduces to the standard Euclidean setting, and the oracle complexity becomes
$\mathcal{O}(d\epsilon^{-3})$, improving the dependence on $\epsilon$ over existing constrained stochastic zeroth-order methods. 
To remove the near-feasibility requirement on the initial point,
we further develop a multi-stage scheme with a rapidly increasing penalty sequence.
%Starting from an arbitrary initial point in $X$, the scheme generates increasingly feasible warm starts and 
It reaches an $\epsilon$-KKT point after
\(\mathcal{O}(1+\log\log(e/\epsilon))\) stages,
and the resulting total oracle complexity preserves the leading high-accuracy complexity of the corresponding single-stage method, with only an additional $\epsilon$-independent first-stage initialization cost.

\end{itemize}

Finally, we report numerical test results on a controlled QCQP dimension-dependence ablation, an MNIST black-box adversarial attack, and a nonconvex fairness-constrained classification problem.
The first two experiments examine the query-efficiency benefit of the non-Euclidean geometry in high dimensions, while the third compares the Euclidean specialization of the proposed method with existing stochastic zeroth-order baselines.

{\bf Organization.} 
In Section \ref{sec:not} notation and preliminaries are introduced. In Section \ref{sec:alg}  a mirror descent linearized augmented Lagrangian method for \eqref{p} is proposed.  In Section \ref{sec4} we present auxiliary lemmas characterizing basic properties of our method. In Section \ref{sec:theo} we conduct a detailed complexity analysis to reach an \(\epsilon\)-KKT point under mirror descent geometries parameterized by $p \ge 2$. In Section \ref{sec:compl} we further provide the specific complexity analysis when using {Rademacher smoothing}.
In Section \ref{sec:num}  numerical results on test problems are reported. Finally, we draw conclusions.

%\newpage
\section{Notation and preliminaries} \label{sec:not}
We use \(\R\) to denote the set of real numbers, \(\R_+\) to denote the set of real numbers greater than or equal to \(0\), and \(\R^d\) to denote the set of \(d\)-dimensional real vectors. Given \(p \ge 2\), we refer to the dual norm of $\|\cdot\|_p$ as $\|\cdot\|_q$, where {$ \frac{1}{p}+ \frac{1}{q}=1$}. For any differentiable function $f$, we denote its gradient by $\nabla f$. 
For iterates generated in the algorithmic process we use the superscript $^k$ to represent the $k$-th iterate and subscript $_i$ to represent $i$-th component of an iterate.
For notation simplicity, given a positive integer $k$ we define $[k] := \{0, 1, \dots, k-1\}$. In expectation subscripts, $u^k:=\{u_j^k\}_j$ and $\xi^k:=\{\xi_j^k\}_j$ denote the collections of random directions and samples drawn at iteration $k$, and we write $u^{[k]}:=\{u^0,\ldots,u^{k-1}\}$ and $\xi^{[k]} := \{\xi^0, \xi^1, \dots, \xi^{k-1}\}$.  Given two sets $A,B \subseteq \R^d$, we define ${\rm dist}_p(A,B) := \inf_{a \in A, b \in B}\|a-b\|_p$ and identify a point with its singleton whenever it appears as an argument of ${\rm dist}_p$. We use \(\mathcal{N}_X\) to denote the normal cone of \(X\) at \(x\), ${\rm I}_d$ to denote the $d$-dimensional identity matrix and $e$ to denots Euler's number.

We next present standard assumptions commonly used in the literature.
Since the zeroth-order oracle queries $F(x+\nu u;\xi)$ even when $x\in X$, the
assumptions on $f$ and $F(\cdot;\xi)$ used below need to hold along the corresponding query segments. For a given 
parameter $\nu>0$ and sampling distribution $\mathbb P_u$, define
\[
    X_\nu:=\{x+t\nu u:\ x\in X,\ t\in[0,1],\ u \sim\mathbb P_u\}.
\]
Let $\mathcal O_\nu$ be an open set containing $X_\nu$.
\begin{assumption}\label{ass:lip}
    The set $X \subset \R^d$ is nonempty, closed, convex and bounded w.r.t. the $q$-norm, i.e., there exists  $D_X>0$ such that $\|x\|_q \le D_X$, for any $x\in X$. The function $f:\R^d \rightarrow \R$ is continuously differentiable over $X$, and  $$\|\nabla f(x)-\nabla f(y)\|_2\le L_f\|x-y\|_2,\quad | f(x) - f(y) | \leq M_f\| x - y\|_2, \quad \forall  x,y\in \mathcal O_\nu,$$
    where $L_f, M_f>0$. 
    Function {$h:\R^d \to \R\cup \{\infty\}$} is proper, lower semicontinuous and convex over $X$, and 
    \[\sup_{g_h \in \partial h(x)} \|g_h\|_2 \leq M_h, \quad  \forall x\in X,\]
    where $M_h>0.$ 
    Functions $c_i:\R^d \rightarrow \R$, $i=1,\ldots,m$ are continuously differentiable over $X$, and 
    $$ |c_i(x)|\le \bar M_c, \quad \|\nabla c_i(x)\|_2\le M_c, \quad \|\nabla c_i(x)-\nabla c_i(y)\|\le L_c\|x-y\|, \quad \forall i=1,\ldots,m;\ x,y\in X,$$
    where $\bar M_c, M_c, L_c>0$. The objective function $\Phi(x)$ has a finite lower bound $\Phi_{\rm inf}$ over $X.$
\end{assumption}

\begin{assumption}\label{ass:u}
  Function \(F(\cdot;\xi)\) is continuously differentiable in
\(x\) for almost every \(\xi\in\Xi\), and 
\[
\E_\xi[\nabla F(x;\xi)]=\nabla f(x), \quad \E_\xi[
\|\nabla F(x;\xi)-\nabla f(x)\|_2^2
]
\le\sigma^2,
\quad
\forall x\in\mathcal O_\nu,
\]
where $\sigma >0.$
\end{assumption}

\begin{assumption}\label{ass:ms-p}
  Function \(F(\cdot;\xi)\) satisfies the mean-squared smoothness on
\(\mathcal O_\nu\):  
\[
(\E_\xi[
\|\nabla F(x;\xi)-\nabla F(y;\xi)\|_2^2
])^{1/2}
\le L_f\|x-y\|_2,
\quad\forall 
x,y\in\mathcal O_\nu.
\]
\end{assumption}

\subsection{Stochastic zeroth-order gradient estimator}\label{ssec:ge}  
The two-point gradient estimator is widely used in  stochastic zeroth-order optimization \cite{duchi2015optimal,sgd,Nesterov17}. Given a parameter $\nu>0$, the two-point gradient estimator of $f$ at $x$ is defined as $\frac{1}{\nu}[f(x+\nu u)-f(x)] u,$
where  $u \in \R^d$ is a random variable following a probability distribution $\mathbb P_u$.
%and $\E_{u\sim \mathbb P} [uu^\top] = {\rm I}_d$.   
However, since the exact objective value of \eqref{p} is not available, we define the stochastic zeroth-order gradient estimator as
\begin{equation}\label{sgrad}
    \begin{aligned}
    G_\nu(x;u,\xi) = \frac{F(x+\nu u;\xi)-F(x;\xi)}{\nu} u, \quad\mbox{where }\ \xi\sim \mathbb P,  u\sim \mathbb P_u.
\end{aligned}\end{equation}
%where the same realization $\xi$ is used at both query points.
For $x\in X$ and $u\sim\mathbb P_u$, both query points and the segment
$\{x+t\nu u:\ t\in[0,1]\}$ lie in $X_\nu\subset\mathcal O_\nu$.
It is easy to check that 
\begin{equation*}
\E_{\xi,u}[G_\nu(x;u,\xi)]=\E_{u\sim \mathbb P_u} [\frac{f(x+\nu u)-f(x)}{\nu} u]=: g_\nu (x).
\end{equation*}
In general, $G_\nu(x;u,\xi)$ is a biased estimator of $\nabla f(x)$. %The random vector $u$ plays a crucial  role in characterizing the property of $G_\nu$.
There are various choices of $u$ following different probability distributions. To unify our  analysis, we give a basic assumption on $u$.

\begin{assumption}\label{ass:exs}
    There exists a positive constant $S_p \ge 1$ such that 
    $\E_{u\sim \mathbb P_u} [\|\langle g , u \rangle u\|_p^2] \leq S_p \|g\|_2^2 $ for any $g \in \R^d
    $. Moreover,  $\E_{u\sim \mathbb P_u} [uu^\top] = {\rm I}_d$ and \(\E_{u\sim \mathbb P_u}[\|u\|_2^6]<\infty\).
\end{assumption}

%\begin{remark}
In the following context, we use abbreviation $\E_{u}$ for $\E_{u\sim \mathbb P_u}$. 
A similar assumption to Assumption \ref{ass:exs} with $p=2$ is required in \cite{duchi2015optimal}.  
In Section \ref{ssec:cpx}, we will present  specific values of $S_p$ corresponding to different probability distributions. Moreover, as $\E_u[u u^\top]={\rm I}_d$, it is straightforward to obtain 
    \begin{equation}\label{le:u1}
    \E_u[\langle g,u\rangle u] = \E_u[uu^\top g] = g,\quad g\in\R^d.
    \end{equation}

We first quantify the bias introduced by the zeroth-order gradient estimator in approximating the true gradient.
\begin{lemma}\label{le:vg}
    Under Assumptions \ref{ass:lip} and  \ref{ass:exs}, it holds that for any $x \in X$, 
    \(
        \|g_\nu(x) - \nabla f(x)\|_p \leq \frac{\nu L_f}{2} \E_u [\|u\|_2^3].
    \)
\end{lemma}
\begin{proof} 
The conclusion follows from the \(L_f\)-smoothness of \(f\), since
        \begin{align*}
           & \|g_\nu(x) - \nabla f(x)\|_p = \|\E_u [\frac{f(x+\nu u)-f(x)}{\nu} u] - \nabla f(x)\|_p \\
            & = \frac{1}{\nu}\|\E_u [(f(x+\nu u)-f(x)- \langle \nabla f(x),\nu u\rangle) u] \|_p\\
            & \leq \frac{1}{\nu} \E_u [|(f(x+\nu u)-f(x)- \langle \nabla f(x),\nu u\rangle)| \| u\|_p]  \leq \frac{\nu L_f}{2} \E_u [\|u\|_2^3],
        \end{align*}
    where the second equality holds from \eqref{le:u1}.
\end{proof}

The next lemma provides mean-square Lipschitz and second-moment bounds for the zeroth-order gradient estimator, which will be used to control the stochastic estimation error.  
\begin{lemma}\label{le:lip}
    Under Assumptions \ref{ass:lip}--\ref{ass:exs}, it holds that for any $x,y \in X$,
    \begin{equation*}     \begin{aligned}
       \E_{u,\xi} [\|G_\nu(y;u,\xi)-G_\nu(x;u,\xi)\|_p^2] & \leq \frac{3\nu^2L_f^2}{2}\E_u [\|u\|_2^6] + 3 S_p L_f^2 \|y-x\|_2^2,\\
         \| g_\nu (y) - g_\nu (x) \|_p^2 & \leq \frac{3\nu^2L_f^2}{2}\E_u [\|u\|_2^6] + 3 S_p L_f^2 \|y-x\|_2^2,\\
        \E_{u,\xi} [\|G_\nu(x;u,\xi)\|_p^2] & \leq \frac{\nu^2L_f^2}{2}\E_u [\|u\|_2^6] + 4 S_p (\sigma^2 + \|\nabla f(x)\|_2^2).
    \end{aligned}
    \end{equation*}
\end{lemma}
\begin{proof}
Note that
    \begin{equation*}     \begin{aligned}
    &\E_{u,\xi} [\|G_\nu(y;u,\xi)-G_\nu(x;u,\xi)\|_p^2]\\
    & =\E_{u,\xi} [\|\frac{F(y+\nu u;\xi)-F(y;\xi)}{\nu} u-\frac{F(x+\nu u;\xi)-F(x;\xi)}{\nu} u\|_p^2]\\
    & =\E_{u,\xi} [\|\frac{F(y+\nu u;\xi)-F(y;\xi)-\langle\nabla F(y;\xi),\nu u\rangle}{\nu} u+\langle \nabla F(y;\xi) - \nabla F(x;\xi), u\rangle u\\
    &\qquad-\frac{F(x+\nu u;\xi)-F(x;\xi)-\langle\nabla F(x;\xi),\nu u\rangle}{\nu} u\|_p^2]\\
    & \leq  \E_u [\frac{3\nu^2L_f^2\|u\|_2^6}{2}] + 3\E_{u,\xi} [\|\langle \nabla F(y;\xi) - \nabla F(x;\xi), u\rangle u\|_p^2]\\
    & \leq  \E_u [\frac{3\nu^2L_f^2\|u\|_2^6}{2}] + 3 S_p \E_\xi [\|\nabla F(y;\xi) - \nabla F(x;\xi)\|_2^2] \\
    &\leq  \E_u [\frac{3\nu^2L_f^2\|u\|_2^6}{2}] + 3 S_p L_f^2 \|y-x\|_2^2,
    \end{aligned}     \end{equation*}
    where the first inequality follows from $(a+b+c)^2\le 3(a^2+b^2+c^2)$ and the following Taylor-remainder estimate: for $z=x,y$, the integral form of the remainder, Minkowski's integral inequality, Cauchy-Schwarz inequality, and Assumption \ref{ass:ms-p} give
    \[
    \E_\xi[\|\frac{F(z+\nu u;\xi)-F(z;\xi)-\langle\nabla F(z;\xi),\nu u\rangle}{\nu}u\|_p^2]
    \leq \frac{\nu^2L_f^2}{4}\|u\|_2^6.
    \]
    The second inequality comes from Assumption \ref{ass:exs}, and the last inequality is derived from Assumption \ref{ass:ms-p}. Secondly, Lemma \ref{le:vg} indicates that
    \begin{align*}
    &\| g_\nu (y) - g_\nu (x) \|_p^2 
    =  \|\E_u [\frac{f(y+\nu u)-f(y)}{\nu} u - \frac{f(x+\nu u)-f(x)}{\nu} u]\|_p^2\\
    & \leq \E_u [\|\frac{f(y+\nu u)-f(y)-\langle \nabla f(y) , \nu u \rangle}{\nu} u + \langle  \nabla f(y) -  \nabla f(x),u\rangle u - \frac{f(x+\nu u)-f(x)-\langle \nabla f(x) , \nu u \rangle}{\nu} u\|_p^2 ]\\
    & \leq   \E_u [\frac{3\nu^2L_f^2\|u\|_2^6}{2}] + 3\E_u [\|\langle \nabla f(y)- \nabla f(x), u\rangle u\|_p^2] \\
    & \leq \E_u [\frac{3\nu^2L_f^2\|u\|_2^6}{2}] + 3S_p L_f^2 \|y-x\|_2^2.
    \end{align*} 
    Thirdly, we also obtain from Assumptions \ref{ass:u}--\ref{ass:exs} that \begin{align*}
            &\E_{u,\xi} [\|G_\nu(x;u,\xi)\|_p^2] = \E_{u,\xi} [|\frac{F(x+\nu u;\xi)-F(x;\xi)}{\nu}|^2 \|u\|_p^2] \\
            & =  \E_{u,\xi} [|\frac{F(x+\nu u;\xi)-F(x;\xi)-\langle \nabla F(x;\xi), \nu u\rangle + \langle \nabla F(x;\xi), \nu u\rangle}{\nu}|^2\|u\|_p^2] \\
             & \le 2\E_{u,\xi} [ \frac{|F(x+\nu u;\xi)-F(x;\xi)-\langle \nabla F(x;\xi), \nu u\rangle |^2}{\nu^2}\|u\|_p^2 + \frac{| \langle \nabla F(x;\xi), \nu u\rangle|^2}{\nu^2}\|u\|_p^2] \\
             & \leq \E_{u,\xi} [\frac{\nu^2L_f^2}{2}\|u\|_2^4\|u\|_p^2 + 2\|\langle \nabla F(x;\xi),u\rangle u \|_p^2]  \leq  \E_{u,\xi} [\frac{\nu^2L_f^2 \| u\|_2^6}{2} + 2S_p\|\nabla F(x;\xi)\|_2^2 ] \\
             & \leq \frac{\nu^2L_f^2}{2}\E_u [\|u\|_2^6] + 4S_p\E_\xi [\|\nabla F(x;\xi)-\nabla f(x)\|_2^2]  + 4S_p\|\nabla f(x)\|_2^2 \\
             & \leq \frac{\nu^2L_f^2}{2}\E_u [\|u\|_2^6] + 4S_p(\sigma^2 + \|\nabla f(x)\|_2^2) .
    \end{align*}
    Here the first term in the second inequality is bounded by the same Taylor-remainder estimate with $z=x$.
    The proof is completed. 
\end{proof}

We next bound the mean-squared  deviation of the stochastic zeroth-order gradient estimator from its expectation, thereby quantifying the stochastic estimation error.
\begin{lemma}\label{le:se}
    Under Assumptions \ref{ass:lip}--\ref{ass:exs}, it holds that 
        \begin{equation*}     \begin{aligned}
            \E_{u,\xi} [\|G_\nu(x;u,\xi)-g_\nu(x)\|_p^2]& \leq 2\nu^2L_f^2\E_u [\|u\|_2^6]+8 S_p (\sigma^2+M_f^2)+ 4M_f^2.
        % \E_{u,\xi} [\|G_\nu(x;u,\xi)-\nabla f(x)\|_p^2]& \leq 5\nu^2L_f^2\E_u [\|u\|_2^6]+16 S_p (\sigma^2+\|\nabla f(x)\|_2^2)+ 8\|\nabla f(x)\|_p^2.
    \end{aligned}     \end{equation*}
\end{lemma}
\begin{proof} From the definition of $G_\nu$ as in \eqref{sgrad}, we  derive
    \bee\begin{aligned}
        \E_{u,\xi} [\|G_\nu(x;u,\xi)-g_\nu(x)\|_p^2] &  \le \E_{u,\xi}[(\|G_\nu(x;u,\xi)\|_p + \|g_\nu(x)-\nabla f(x)\|_p + \|\nabla f(x)\|_p)^2] \\
        &\leq \E_{u,\xi} [2\|G_\nu(x;u,\xi)\|_p^2 +4\|g_\nu(x)-\nabla f(x)\|_p^2 + 4\|\nabla f(x)\|_p^2]\\
        & \leq 2\nu^2L_f^2\E_u [\|u\|_2^6]+8 S_p (\sigma^2+\|\nabla f(x)\|_2^2)+ 4\|\nabla f(x)\|_p^2,
    \end{aligned}  \eee
   where the second inequality uses $(a+b+c)^2\le 2a^2 + 2(b+c)^2\le 2 a^2 + 4b^2 + 4c^2$ and  the last inequality uses  Lemmas \ref{le:vg}, \ref{le:lip}. Thus, we derive the estimate using \(\|\cdot\|_p \le \|\cdot\|_2\) together with the fact that the \(M_f\)-Lipschitz continuity of  differentiable \(f\) on the open set \(\mathcal O_\nu\) implies \(\|\nabla f(x)\|_2 \le M_f\).
\end{proof}

\section{Mirror descent linearized augmented Lagrangian method}\label{sec:alg}

The augmented Lagrangian (AL) function associated with \eqref{p} is given by 
\begin{equation*}
\begin{aligned}
    \mathcal{L}(x,\lambda;\mu) = f(x) +  h(x) + \lambda^\top c(x) + \frac{\mu}{2} \|c(x)\|_2^2,
\end{aligned} \end{equation*}
where 
$\mu>0$ is a penalty parameter and $\lambda \in \R^m$. 
For notation brevity, we define 
\begin{align*}
\phi(x,\lambda;\mu)&=f(x) + \lambda^\top c(x) + \frac{\mu}{2} \|c(x)\|_2^2.
\end{align*} 
At current iterate \(x^k\), instead of minimizing the AL function to update the primal variable \(x\), we construct a simpler approximation model motivated by the mirror descent mechanism:
\[ \mathcal L(x,\lambda; \mu)\approx \phi(x^k,\lambda;\mu) + \langle \nabla_x \phi(x^k,\lambda;\mu),x \rangle +  h(x) + \frac{1}{\eta}V(x^k,x),
\]
where $V(x,y) := v(y) - v(x) -\langle \nabla v(x) , y-x \rangle$ for $x,y\in \R^d$ is the Bregman divergence specifying the mirror descent step. Here, we impose the following assumption on the distance-generating function $v$. An available choice of $v$ will be given in Section \ref{ssec:cpx}.
\begin{assumption}\label{ass:v}
    The distance-generating function $v: \mathbb{R}^d \to \mathbb{R}$ is $1$-strongly convex   w.r.t. the \(\ell_q\)-norm, meaning 
\be\label{v_str_cvx}
\langle x - y, \nabla v(x) - \nabla v(y) \rangle \geq \|x - y\|_q^2, \quad  x,y\in \mathbb R^d.
\ee
And $v$ is $L_v$-smooth over $X$ w.r.t. the $\ell_q$-norm, that is
\be\label{v_Lv_sm}
\|\nabla v(x)-\nabla v(y)\|_p
    \le L_v\|x-y\|_q,\quad  x,y \in X. 
\ee
\end{assumption}
% When smoothness of \(v(x)\) is required, $L_v$-smoothness within the constraint set \(X\) w.r.t. \ the \(\ell_q\)-norm means
% \[
%     \|\nabla v(x)-\nabla v(y)\|_p
%     \le L_v\|x-y\|_q,
%     \qquad x,y\in X.
% \]
Note that  $\nabla_x\phi(x^k,\lambda; \mu)$ cannot be obtained exactly as the exact gradient of $f$ is not available. We therefore compute a stochastic gradient estimate of $f$ at $x^k$, defined as
\begin{align}\label{hsgd}
            s^k = \begin{cases}
            \frac{1}{n_0}\sum_{j=1}^{n_0}G_\nu(x^0;u_j^0,\xi_j^0), & k=0,\\
            \frac{1}{n}\sum_{j=1}^n [G_\nu(x^k;u_j^k,\xi_j^k) + (1-\alpha)( s^{k-1} - G_\nu(x^{k-1};u_j^k,\xi_j^k))], \quad & k\ge 1,
            \end{cases}
\end{align}
where  
%conditioned on the history before iteration $k$, the pairs of random directions and random samples $\{(u_j^k,\xi_j^k)\}_{j=1}^n$ are mutually independent for $k\ge1$, and $\{(u_j^0,\xi_j^0)\}_{j=1}^{n_0}$ are mutually independent for $k=0$. In particular, 
the random directions $\{u_j^k\}_{j,k}\sim\mathbb P_u$ and random samples $\{\xi_j^k\}_{j,k}\sim \mathbb P$ are mutually  independent and independent of $x$. To reduce the variance caused by the stochasticity of the objective and zeroth-order approximations (see Lemma \ref{le:se}), we incorporate the momentum technique  \cite{storm,tran2022hybrid}. Meanwhile, the zeroth-order approximation introduces a dimension-dependent second moment and also affects the mean-square Lipschitz bound of the approximate gradient (see Lemma \ref{le:lip}). We therefore use a mini-batch sampling whose size reflects this dimension dependence.

The complete algorithm framework is summarized in Algorithm \ref{algo}.

\begin{algorithm}[ht]
    \caption{Mirror  descent linearized augmented Lagrangian method}
    \label{algo}
    \begin{algorithmic}[1]
    \REQUIRE{$x^0\in X, \lambda^0 = 0,$ positive parameters $ \mu, \eta, \nu, \alpha \in (0,1] , $ $\{\rho_k\in(0,\mu)\}$} and positive integers  $K, n_0, n$. 
    \ENSURE $x^{R}$, where $R \in [K]$ is uniformly chosen at random.
    \FOR{$k=0,\ldots,K-1$}
    \STATE{Compute $s^k$ through \eqref{hsgd}.
    }
    \STATE{Compute $x^{k+1}$ through
        \begin{align}\label{x}
            x^{k+1} = \arg\min_{x \in X} \, \langle s^k +\nabla c(x^k)^\top(\lambda^k+\mu c(x^k)),x \rangle +  h(x) + \frac{1}{\eta}V(x^k,x).
        \end{align}
    }
    \STATE{Compute $\lambda^{k+1}$ through
        $%\label{multiplier}
            \lambda^{k+1} = \lambda^k + \rho_k c(x^k).
        $
    }
    \ENDFOR
    \end{algorithmic}
\end{algorithm}

In order to measure the stationarity of the output  returned by Algorithm \ref{algo}, we employ the generalized gradient mapping, defined by
\begin{equation*}     \begin{aligned}
    \mathcal{G}_X (x,g,\eta) = \frac{1}{\eta}(x-x^+), \quad \mbox{where $x^+ = \arg \min_{y \in X} \langle g,y \rangle + h(y) + \frac{1}{\eta} V(x,y).$}
\end{aligned}  \end{equation*}  
The following lemma shows that the generalized gradient mapping is consistent with the standard  KKT stationarity condition.

\begin{lemma}
    Given $x \in \R^d$ and $\lambda \in \R^m$, if $\mathcal{G}_X (x,\nabla f(x) + \nabla c(x)^\top \lambda,\eta) = 0$, it holds that  
    \bee
        0 \in \nabla f(x) + \nabla c(x)^\top \lambda + \partial  h(x) + \mathcal{N}_X(x).
    \eee
\end{lemma}
\begin{proof}
    Let 
    \(
    x^+ = \arg \min_{y \in X} \langle \nabla f(x)+ \nabla c(x)^\top\lambda,y \rangle + h(y) + \frac{1}{\eta} V(x,y).
    \)
    Then from $x^+ = x$ and  the first-order optimality condition at $x^+$, 
    we obtain the conclusion. 
\end{proof}

Our goal in this paper is to establish the oracle complexity of Algorithm \ref{algo} for returning a random $\epsilon$-KKT point of \eqref{p}, which is defined as follows.
\begin{definition}{\label{de:kkt}}
    Given $\epsilon > 0$, a point $x\in X$ is called an
    \(\epsilon\)-KKT point of \eqref{p}, whose stationarity and feasibility are
    measured in the $\ell_q$- and $\ell_p$-norms, respectively, if exists $\lambda\in\mathbb{R}^m$ such that 
    \[
         \mathbb{E} [\|\mathcal{G}_X ( x,\nabla f(x) + \nabla c(x)^\top\lambda,\eta)\|_q^2  ] \le \epsilon^2, \quad   
        \mathbb{E} [\|c(x)\|_p^2]  \le \epsilon^2,
    \]
    where the expectation is taken w.r.t.  all random variables generated by the algorithm.
\end{definition}

To establish a near-feasibility guarantee in expectation for the output of Algorithm \ref{algo}, we impose the following regularity assumption on its iterates.  

\begin{assumption}\label{ass:ensc}
   % (Extended NonSingularity Condition) 
   There exists a constant $ \beta > 0 $ such that
    \[
            \beta \| c ( x^k  )   \|_p \leq  {\rm dist}_p ( {\nabla c} ( x^k )^\top c ( x^k ),  -\mathcal{N}_{X}(x^{k})  )  ,\quad \forall k \in  [ K+1  ] .
    \]
\end{assumption}

Assumption \ref{ass:ensc} is an extension of NonSingularity Condition (NSC) \cite[Assumption 4]{ippp} corresponding to \(p=2\). NSC and its variants have been used to analyze  complexities of AL methods in recent years \cite{rialm2021,li2024stochastic,ippp,sahin2019inexact}. 
Such a regularity condition provides a quantitative relation between the
constraint violation and the corresponding constraint-gradient residual,
which plays a key role in establishing near-feasibility from approximate
stationarity of penalized or AL subproblems
\cite{li2024stochastic,ippp}.
%This type of regularity condition is conceptually different from classical constraint qualifications such as the strong MFCQ; in general, neither condition is strictly stronger nor strictly weaker than the other \cite{li2024stochastic,ippp}.
%In contrast to the global regularity condition imposed over the constraint set in \cite{li2024stochastic}, Assumption~\ref{ass:ensc} only requires the regularity condition to hold uniformly along the iterates generated by Algorithm~\ref{algo}.

When $p=2$, with the Euclidean distance-generating
function $v(x)=\frac{1}{2}\|x\|_2^2$, Definition~\ref{de:kkt} reduces to the
widely used $\epsilon$-KKT definition in  $\ell_2$-norm,  defined by 
\be\label{eps-KKT-2}
\mathbb{E} [\|\mathcal{G}_X (x,\nabla f(x) + \nabla c(x)^\top\lambda,\eta)\|_2^2  ] \le \epsilon^2, \mbox{ and } \E[\|c(x)\|_2^2]\le \epsilon^2.
\ee
and the regularity assumption  characterized through 
 \begin{align}
\exists \beta>0\  \mbox{ s.t. }\  \beta \| c ( x^k  )   \|_2 \leq  {\rm dist}_2 ( {\nabla c} ( x^k )^\top c ( x^k ),  -\mathcal{N}_{X}(x^{k})  )  ,\quad \forall k \in  [ K+1  ] .\label{CQ_2}
 \end{align}
We will show in
Corollary~\ref{thm:2-norm} the corresponding oracle complexity of the proposed algorithm under this Euclidean specialization.

\section{Auxiliary lemmas} \label{sec4}
In this section, we present  auxiliary lemmas that are used in subsequent sections.  We first recall a smoothness inequality for the squared $\ell_p$-norm, which will be used to control sums of stochastic estimation errors.
\begin{lemma}\label{le:normp}
    For any $p \geq 2$,  it holds that for all $x,y \in \R^d$,
    \begin{align*}
        \|x+y\|_p^2 \leq \|x\|_p^2 + \langle g_x ,y \rangle + (p-1) \|y\|_p^2,\quad\mbox{where $g_x = \nabla \|\cdot\|_p^2(x)$.}
    \end{align*}  
\end{lemma}
\begin{proof}
    The result comes from the fact that function $\frac{1}{2}\|x\|_p^2$ is $(p-1)$-smooth w.r.t. the \(\ell_p\)-norm \cite{beck2017first}.
\end{proof}

We establish uniform bounds on the multiplier sequence generated by Algorithm \ref{algo} and its variation.
\begin{lemma}\label{new:z2}
    Under Assumption \ref{ass:lip}, it holds that  
    \[|\lambda_i^k| \le \bar{M}_c \sum_{t=0}^{k-1} \rho_t,\quad|\lambda_i^{k} - \lambda_i^{k-1}| \le \rho_{k-1} \bar{M}_c,\quad \forall i=1,\ldots,m;\ k = 1, \dots, K.\] 
      
\end{lemma}
\begin{proof}
    The result is due to \(\lambda^{k+1} = \lambda^k + \rho_k c(x^k), k=0,\ldots,K-1\) and \(|c_i(x)| \leq \bar{M}_c, \forall i=1,\ldots,m\), together with \(\lambda_i^0 = 0\).
\end{proof}

Using the multiplier bounds above, the next lemma establishes the Lipschitz continuity of $\nabla_x\phi(x,\lambda^k;\mu)$.
\begin{lemma}\label{le:Lu}
    Under Assumptions \ref{ass:lip}--\ref{ass:u},  it holds that for any $x,y \in X $, $\mu>0$ and $k =0,\ldots,K$, 
    \bee
        \lVert \nabla_x {\phi} ( x, \lambda^k;\mu ) - \nabla_x {\phi} ( y, \lambda^k; \mu)  \rVert_2
        \leq L_\mu \lVert x-y \rVert_2 ,
    \eee
where \(L_\mu = L_f + m(\mu M_c^2+\mu \bar{M}_c L_c + \rho \bar{M}_c L_c) \) with \(\rho\geq\sum_{k=0}^{K-1} \rho_k\).
\end{lemma}
\begin{proof}
     It follows from Assumptions \ref{ass:lip}--\ref{ass:u} that
    \bee
    \begin{aligned}
        &\lVert \nabla_x {\phi} ( x, \lambda^k;\mu ) - \nabla_x {\phi} ( y, \lambda^k; \mu)  \rVert_2 \\
        &\leq  \lVert \nabla f(x) - \nabla f(y) \rVert_2 + \sum_{i = 1}^m  \lVert  ( \mu  c_{i} ( x ) +\lambda_i^k ) \nabla c_{i} ( x ) - ( \mu  c_{i} ( y ) +\lambda_i^k )\nabla c_{i} ( y )   \rVert_2 \\
        &\leq   {L_f} \|x-y\|_2+ \sum_{i = 1}^m  \lVert (   \mu  c_{i} ( x )  -  \mu  c_{i} ( y ) ) \nabla c_{i} ( x ) + ( \mu  c_{i} ( y ) +\lambda_i^k )  ( \nabla c_i ( x ) -\nabla c_i ( y )  )   \rVert_2 \\
        &\leq   {L_f} \|x-y\|_2+ \sum_{i = 1}^m  (  \mu \lvert c_{i} ( x ) -c_{i} ( y )  \rvert \lVert \nabla c_{i} ( x )  \rVert_2+ | \mu  c_{i} ( y ) +\lambda_{i}^k | L_c\lVert x-y\rVert_2  )\\
        &\leq   ({L_f} + m(\mu M_c^2+\mu \bar{M}_c L_c + \rho \bar M_c L_c)) \|x-y\|_2,
    \end{aligned}
    \eee
    which derives the conclusion.
   %The last inequality follows from \(\|\nabla c_i(x)\|_2\le M_c\); this follows from the boundedness of \(X\) and the \(L_c\)-Lipschitz continuity of \(\nabla c_i\), where \(M_c\) may be enlarged if necessary. 
\end{proof}
The following lemma establishes a one-step descent inequality for the AL function.
\begin{lemma}\label{le:x} 
     Under Assumptions \ref{ass:lip}--\ref{ass:exs} and Assumption \ref{ass:v}, 
     %\red{suppose \(v(x)\) is \(1\)-strongly convex w.r.t. \(\ell_q\)-norm.} 
     it holds that for any \(k=0,\ldots,K-1\),
    \bee
    \begin{aligned}
    \frac{1}{2\eta}\|x^{k+1}-x^k\|_q^2-\frac{L_\mu}{2}\|x^{k+1}-x^k\|_2^2 &\leq \mathcal{L}(x^k,\lambda^k;\mu)-\mathcal{L}(x^{k+1},\lambda^{k+1};\mu) + \frac{\eta \nu^2 L_f^2}{4}\E_u[\|u\|_2^6] \\
    &\quad + \eta \|s^k - g_\nu (x^k)\|_p^2 + m \rho_k \bar{M}_c^2.
    \end{aligned}
    \eee
\end{lemma}
\begin{proof}
    It follows from Assumption \ref{ass:lip} and Lemma \ref{new:z2} that
    \bee
    \begin{aligned}
        &\mathcal{L} ( x^{k+1},\lambda^{k};\mu )  = \mathcal{L} ( x^{k+1},\lambda^{k+1} ;\mu) +\sum_{i=1}^m   (\lambda_i^k -  \lambda_i^{k+1}  )  c_i ( x^{k+1} )    \\
        &\geq  \mathcal{L} ( x^{k+1},\lambda^{k+1} ;\mu) -\sum_{i=1}^m   | \lambda_i^{k}-\lambda_i^{k+1}   | | c_i ( x^{k+1} )   |
         \ge \mathcal{L} ( x^{k+1},\lambda^{k+1};\mu ) -m \rho_k \bar{M}_c^2.
    \end{aligned}
    \eee
    By the optimality condition for \eqref{x}, there exists  $g_h^{k+1} \in \partial h(x^{k+1})$ such that
    \bee
         \langle s^k+\nabla c(x^k)^\top (\lambda^k+\mu c(x^k))+g_h^{k+1}+\frac{1}{\eta}(\nabla v(x^{k+1})-\nabla v(x^k)),x-x^{k+1} \rangle \geq 0,\quad \forall x \in X.
    \eee
    Then by  the convexity of $h $ and the setting $x=x^k $, we have
    \bee
    \begin{aligned}
        &h(x^{k+1})-h(x^k) \leq \langle g_h^{k+1}, x^{k+1}-x^k \rangle \\
        &\leq - \langle s^k+\nabla c(x^k)^\top (\lambda^k+\mu c(x^k))+\frac{1}{\eta}(\nabla v(x^{k+1})-\nabla v(x^k)),x^{k+1}-x^k \rangle,
    \end{aligned}
    \eee
     which, together with Lemma \ref{le:Lu}, implies that
    \begin{align*}
        &\mathcal{L}(x^{k+1},\lambda^{k+1};\mu)-\mathcal{L}(x^k,\lambda^k;\mu)\\
        & = \mathcal{L}(x^{k+1},\lambda^{k+1};\mu) -\mathcal{L}(x^{k+1},\lambda^k;\mu)+\mathcal{L}(x^{k+1},\lambda^k;\mu)-\mathcal{L}(x^k,\lambda^k;\mu)\\
        & \leq \mathcal{L}(x^{k+1},\lambda^k;\mu) -\mathcal{L}(x^k,\lambda^k;\mu)+m\rho_k \bar{M}_c^2\\
        & \leq  \langle \nabla_x \phi(x^k,\lambda^k;\mu),x^{k+1}-x^k\rangle + \frac{L_\mu}{2}\|x^{k+1}-x^k\|_2^2+h(x^{k+1})-h(x^k)+m\rho_k \bar{M}_c^2\\
        & \leq \langle \nabla f (x^k),x^{k+1}-x^k\rangle - \langle s^k+\frac{1}{\eta}(\nabla v(x^{k+1})-\nabla v(x^k)),x^{k+1}-x^k\rangle + \frac{L_\mu}{2}\|x^{k+1}-x^k\|_2^2+m\rho_k \bar{M}_c^2\\
        &= \langle \nabla f(x^k)-g_\nu(x^k)+g_\nu (x^k) - s^k - \frac{1}{\eta}(\nabla v(x^{k+1})-\nabla v(x^k)),x^{k+1}-x^k\rangle  + \frac{L_\mu}{2}\|x^{k+1}-x^k\|_2^2+m\rho_k \bar{M}_c^2\\
        & \leq \eta \|\nabla f(x^k)-g_\nu(x^k)\|_p^2 + \eta \|s^k - g_\nu (x^k)\|_p^2-\frac{1}{2\eta} \|x^{k+1}-x^k\|_q^2 + \frac{L_\mu}{2}\|x^{k+1}-x^k\|_2^2 +m\rho_k \bar{M}_c^2\\
        & \leq \frac{\eta \nu^2 L_f^2}{4}\E_u[\|u\|_2^6] + \eta \|s^k - g_\nu (x^k)\|_p^2 -\frac{1}{2\eta} \|x^{k+1}-x^k\|_q^2 + \frac{L_\mu}{2}\|x^{k+1}-x^k\|_2^2+m\rho_k \bar{M}_c^2,
    \end{align*}
    where the last line uses Lemma \ref{le:vg}. By arranging the terms we obtain the conclusion. 
\end{proof}

To better understand the behavior of $s^k$, we define the error 
\bee
\varepsilon^k := s^k - g_\nu (x^k), \quad k=0,\ldots,K-1.
\eee  
and
    \be \label{ab}
    \begin{aligned}
        \tilde \sigma^2 &= \nu^2 L_f^2 \E_u[\|u\|_2^6]+4 S_p (\sigma^2+M_f^2)+ 2M_f^2,\quad
          \Delta_0 = \mathcal{L}(x^0,\lambda^0;\mu)-\Phi_{\rm inf}+2m \rho \bar{M}_c^2.
    \end{aligned}
    \ee
   %  where $\Phi_{\rm inf}$ is a lower bound for $f$ over $X$. 

  The following lemma derives a recursive bound for the mean-square error of the momentum-based zeroth-order gradient estimator.
\begin{lemma}\label{le:eps}
    Under Assumptions \ref{ass:lip}--\ref{ass:exs}, let \(K\ge2\). Then it holds that  for $k=0,\ldots,K-2$,\ %with $p \in [2,{2\ln d}]$, 
   
    \be\label{ave-eps}
    \begin{aligned}
         &\E_{u^{[k+2]},\xi^{[k+2]}}[\|\varepsilon^{k+1}\|_p^2]
    \leq (1-\alpha)^2 \E_{u^{[k+1]},\xi^{[k+1]}}[\|\varepsilon^k\|_p^2]\\
    &+\frac{p-1}{n}[4\alpha^2 {\tilde \sigma^2} + 12 (1-\alpha)^2 \nu^2 L_f^2 \E_u[\|u\|_2^6] + 24(1-\alpha)^2S_pL_f^2
\E_{u^{[k+1]},\xi^{[k+1]}}
[\|x^{k+1}-x^k\|_2^2]].
    \end{aligned}
    \ee
    %where $\tilde \sigma^2$ is introduced in \eqref{ab}.
\end{lemma}
\begin{proof}
Let $U^{k+1}:=(u_1^{k+1},\ldots,u_n^{k+1})$ and $\Xi^{k+1}:=(\xi_1^{k+1},\ldots,\xi_n^{k+1})$. 
    First, by the definition of $\varepsilon^k$, we have
    \begin{align*}
        \varepsilon^{k+1} &= s^{k+1} - g_\nu (x^{k+1})
             = \frac{1}{n}\sum_{j=1}^n [ G_\nu (x^{k+1};u_j^{k+1},\xi_j^{k+1}) - g_\nu (x^{k+1}) + (1-\alpha) (s^k - G_\nu (x^k;u_j^{k+1},\xi_j^{k+1}))]\\
            &= \frac{1}{n}\sum_{j=1}^n [ G_\nu (x^{k+1};u_j^{k+1},\xi_j^{k+1}) - g_\nu (x^{k+1})+ (1-\alpha) (\varepsilon^k +g_\nu (x^k) - G_\nu (x^k;u_j^{k+1},\xi_j^{k+1}))]. 
    \end{align*} 
    Conditioning on $u^{[k+1]}$ and $\xi^{[k+1]}$, and taking 
expectation w.r.t. the new batch
$(U^{k+1},\Xi^{k+1})$ yield
    \begin{align}
        \E_{U^{k+1},\Xi^{k+1}}[\|\varepsilon^{k+1}\|_p^2]
        &= \E_{U^{k+1},\Xi^{k+1}} [\|(1-\alpha) \varepsilon^k+\frac{1}{n}\sum_{j=1}^n G_\nu (x^{k+1};u_j^{k+1},\xi_j^{k+1}) - g_\nu (x^{k+1})\notag  \\
        &\quad + (1-\alpha) (g_\nu (x^k) - \frac{1}{n}\sum_{j=1}^n G_\nu (x^k;u_j^{k+1},\xi_j^{k+1}))\|_p^2]\notag  \\
        & \leq (1-\alpha)^2 \|\varepsilon^k\|_p^2+\frac{p -1}{n^2} \E_{U^{k+1},\Xi^{k+1}} [\sum_{j=1}^n\| G_\nu (x^{k+1};u_j^{k+1},\xi_j^{k+1}) - g_\nu (x^{k+1}) \notag \\
        &\quad+ (1-\alpha) (g_\nu (x^k) - G_\nu (x^k;u_j^{k+1},\xi_j^{k+1}))\|_p^2], \label{ave-eps-k}
    \end{align}
   where the inequality follows from Lemma \ref{le:normp} with $x=(1-\alpha)\varepsilon^k$, $y = \varepsilon^{k+1} -x$ and $\E_{U^{k+1},\Xi^{k+1}}[y]={0}$. We now focus on the second term of the right side of \eqref{ave-eps-k}. Note that 
    \begin{align*}
        &\E_{u_j^{k+1},\xi_j^{k+1}} [\|G_\nu (x^{k+1};u_j^{k+1},\xi_j^{k+1}) - g_\nu (x^{k+1}) + (1-\alpha) (g_\nu (x^k) - G_\nu (x^k;u_j^{k+1},\xi_j^{k+1}))\|_p^2]\\
        &\leq  \E_{u_j^{k+1},\xi_j^{k+1}} [2\alpha^2\|G_\nu (x^{k+1};u_j^{k+1},\xi_j^{k+1}) - g_\nu (x^{k+1})\|_p^2\\
        &\quad + 2(1-\alpha)^2\|G_\nu (x^{k+1};u_j^{k+1},\xi_j^{k+1}) - g_\nu (x^{k+1})+g_\nu (x^k) - G_\nu (x^k;u_j^{k+1},\xi_j^{k+1})\|_p^2]\\
        & \leq  \E_{u_j^{k+1},\xi_j^{k+1}} [2\alpha^2\|G_\nu (x^{k+1};u_j^{k+1},\xi_j^{k+1}) - g_\nu (x^{k+1})\|_p^2 + 4(1-\alpha)^2\|g_\nu (x^k) - g_\nu (x^{k+1}) \|_p^2\\
        & \quad + 4(1-\alpha)^2\|G_\nu (x^{k+1};u_j^{k+1},\xi_j^{k+1}) - G_\nu (x^k;u_j^{k+1},\xi_j^{k+1})\|_p^2]\\
        & \leq  4\alpha^2 (\nu^2L_f^2 \E_u[\|u\|_2^6] +4 S_p (\sigma^2+M_f^2)+ 2M_f^2)  + 12 (1-\alpha)^2 \nu^2 L_f^2 \E_u[\|u\|_2^6] + 24 (1-\alpha)^2  S_p  L_f^2  \|x^{k+1}-x^k\|_2^2\\
        & =  4\alpha^2 {\tilde \sigma^2} + 12 (1-\alpha)^2 \nu^2 L_f^2 \E_u[\|u\|_2^6] + 24 (1-\alpha)^2  S_p  L_f^2  \|x^{k+1}-x^k\|_2^2,
    \end{align*} 
    where the first and second inequalities come from $\|a+b\|_p \leq \|a\|_p + \|b\|_p$ and $2\|a\|_p\|b\|_p \leq \|a\|_p^2+\|b\|_p^2$, the third one follows from Lemma \ref{le:lip} and Lemma \ref{le:se}, and  the equality is due to the definition of \(\tilde \sigma^2\). Taking expectation w.r.t. 
$u^{[k+1]}$ and $\xi^{[k+1]}$ on both sides and using
the tower property yields \eqref{ave-eps}. 
\end{proof}

\section{Oracle complexity analysis} \label{sec:theo}
In this section, we establish the oracle complexity of Algorithm \ref{algo}, in terms of the total number of stochastic function evaluations, to return an $\epsilon$-KKT point of \eqref{p}. In Subsection \ref{sec3-1}, we present upper bound estimates on the optimality measure in expectation. In Subsection~\ref{sec3-2}, we derive oracle complexity bounds for
the two cases corresponding to \(p\in[2,2\ln d]\) and \(p>2\ln d\), respectively. In Subsection \ref{ssec:cpx}, we provide discussions on the choice of the mirror descent geometry and the values of \(S_p\).

\subsection{Upper bound estimates on KKT measure} \label{sec3-1}

\begin{theorem}[{\bf {Stationarity}}]\label{stationary}
     Under Assumptions \ref{ass:lip}--\ref{ass:exs} and Assumption \ref{ass:v}, let \(K\ge2\) and 
    \bee
     \rho_k = \frac{\rho}{K},\ \eta L_\mu \leq \frac{1}{2},\ 0<\alpha\le1~{\rm and}~\ n \ge \frac{96(p-1)S_p\eta^2L_f^2}{\alpha} 
    \eee
    for $k=0,\ldots, K-1$. 
    Then it holds that with \(\tilde \lambda^R = \lambda^R + \mu c(x^R)\),
    \be \label{bound}
    \begin{aligned}
        &\E_{R;u^{[K]},\xi^{[K]}}  [\|\mathcal{G}_X (x^R,\nabla f(x^R) + \nabla c(x^R)^\top\tilde \lambda^R,\eta)\|_q^2 ] \\
        & \leq (4L_f^2+3\eta^{-2})\nu^2 \E_u[\|u\|_2^6] +\frac{30(p-1)\tilde \sigma^2}{n_0 \alpha K}+\frac{60\alpha (p -1)\tilde \sigma^2}{n} +\frac{27 \Delta_0}{ \eta K}.
    \end{aligned}
    \ee 
   % where $\tilde \sigma^2$ and $\Delta_0$ are introduced in \eqref{ab}.
\end{theorem}
\begin{proof}
    First,  {due to $\|\cdot\|_2\leq\|\cdot\|_q$, it follows from} $\eta L_\mu \leq \frac{1}{2}$ that 
    \[
    \frac{1}{2\eta}\|x^{k+1}-x^k\|_q^2 - \frac{L_\mu}{2}\|x^{k+1}-x^k\|_2^2\ge (\frac{1}{2\eta}-\frac{L_\mu}{2})\|x^{k+1}-x^k\|_q^2 \ge \frac{1}{4\eta}\|x^{k+1}-x^k\|_q^2,
    \]
    which further yields from Lemma \ref{le:x} that 
    \bee 
    \begin{aligned}
        \|x^{k+1}-x^k\|_q^2 &\leq  4\eta ( \mathcal{L}(x^k,\lambda^k;\mu)-\mathcal{L}(x^{k+1},\lambda^{k+1};\mu)+ m  \rho_k \bar{M}_c^2) + \eta^2 \nu^2 L_f^2\E_u[\|u\|_2^6] +4\eta^2\|\varepsilon^k\|_p^2.
    \end{aligned}
    \eee
    Taking expectation w.r.t. \(u^{[k+1]},\xi^{[k+1]}\) on both sides of above inequality and 
    summing the above inequality over $k=0,\ldots, K-1$ yields
    \begin{align*}
     \sum_{k=0}^{K-1} \E_{u^{[k+1]},\xi^{[k+1]}} \|x^{k+1}-x^k\|_q^2 &\leq 4\eta (\mathcal L(x^0,\lambda^0;\mu)-\E_{u^{[k+1]},\xi^{[k+1]}}[\mathcal L(x^{k},\lambda^k;\mu)])+4K\eta m\sum_{k=0}^{K-1}\rho_k \bar{M}_c^2
     \\&\quad +K\eta^2 \nu^2 L_f^2\E_u[\|u\|_2^6] +4K\eta^2\sum_{k=0}^{K-1} \E_{u^{[k+1]},\xi^{[k+1]}}[\|\varepsilon^k\|_p^2].
   \end{align*}
    Note that for any \(t=1,\ldots,K\), 
    \be\label{del}
    \begin{aligned}
    \Delta_0 =\mathcal{L}(x^0,\lambda^0;\mu) - \Phi_{\rm inf}+2m \rho \bar{M}_c^2 \geq \mathcal{L}(x^0,\lambda^0;\mu)
    -\E_{u^{[t]},\xi^{[t]}}[\mathcal{L}(x^t,\lambda^t;\mu)]
    +m \sum_{k=0}^{t-1}\rho_k \bar{M}_c^2,
    \end{aligned}
    \ee
    where the inequality follows from $\Phi(x)\ge\Phi_{\inf}$, Lemma \ref{new:z2}, $|c_i(x^t)|\le \bar{M}_c$, and $\sum_{k=0}^{t-1}\rho_k\le\rho$, which give $\mathcal{L}(x^t,\lambda^t;\mu)=\Phi(x^t)+(\lambda^t)^\top c(x^t)
+\frac{\mu}{2}\|c(x^t)\|_2^2\geq\Phi_{\inf}-m\rho \bar{M}_c^2$. Therefore, we obtain
    \be\label{edx}
    \begin{aligned}
        \frac{1}{K} \sum_{k=0}^{K-1} \E_{u^{[k+1]},\xi^{[k+1]}} \|x^{k+1}-x^k\|_q^2 
        \leq \frac{4\eta \Delta_0}{K} + \eta^2 \nu^2 L_f^2 \E_u[\|u\|_2^6] +\frac{4 \eta^2 }{K} \sum_{k=0}^{K-1} \E_{u^{[k+1]},\xi^{[k+1]}}[\|\varepsilon^k\|_p^2],
    \end{aligned}
    \ee
       It then, together with $\|\cdot\|_2\le \|\cdot\|_q$ and Lemma \ref{le:eps}, yields
    \be\label{edx2}
    \begin{aligned}
        &\E_{u^{[k+2]},\xi^{[k+2]}}[\|\varepsilon^{k+1}\|_p^2] 
         \leq (1-\alpha)^2 \E_{u^{[k+1]},\xi^{[k+1]}}[\|\varepsilon^k\|_p^2] +\frac{4(p-1)}{n}[\alpha^2\tilde \sigma^2+3\nu^2L_f^2\E_u[\|u\|_2^6]]\\
        &\qquad\qquad\qquad\qquad\qquad\qquad + \frac{24(p-1) (1-\alpha)^2 S_p L_f^2}{n} \E_{u^{[k+1]},\xi^{[k+1]}}[\|x^{k+1}-x^k\|_2^2]\\
        & \leq (1-\alpha)^2(1+\frac{96(p -1)\eta^2 S_p L_f^2 }{n}) \E_{u^{[k+1]},\xi^{[k+1]}}[\|\varepsilon^k\|_p^2] +4\alpha^2  (p -1) \frac{\tilde \sigma^2}{n} + \frac{12(p-1)\nu^2L_f^2\E_u[\|u\|_2^6]}{n}(1+\frac{S_p}{2})\\
        &\quad + \frac{96 (p -1) \eta S_p L_f^2}{n} (\E_{u^{[k+1]},\xi^{[k+1]}}[\mathcal{L}(x^k,\lambda^k;\mu)-\mathcal{L}(x^{k+1},\lambda^{k+1};\mu)]+m \rho_k \bar{M}_c^2)\\
        &\leq  (1-\alpha) \E_{u^{[k+1]},\xi^{[k+1]}}[\|\varepsilon^k\|_p^2]+4\alpha^2  (p -1)\frac{\tilde \sigma^2}{n} + \frac{18(p-1)S_p\nu^2L_f^2\E_u[\|u\|_2^6]}{n}\\ 
        &\quad + \frac{96 (p -1) \eta S_p L_f^2}{n} (\E_{u^{[k+1]},\xi^{[k+1]}}[\mathcal{L}(x^k,\lambda^k;\mu)-\mathcal{L}(x^{k+1},\lambda^{k+1};\mu)]+m \rho_k \bar{M}_c^2), 
    \end{aligned}
    \ee
    where the first term in the last inequality follows from \(0<\alpha\le1\) and {$n \ge \frac{96(p-1)S_p\eta^2L_f^2}{\alpha}$}, and the third term is due to \(S_p \ge 1\). Let \(Z_j^0:=G_\nu (x^0;u_j^0,\xi_j^0)-g_\nu(x^0)\). Applying Lemma \ref{le:normp} successively to the partial sums and conditioning on \(\{Z_i^0\}_{i<j}\), the linear term has zero expectation as  the vectors \(\{Z_j^0\}_{j=1}^{n_0}\) are independent and mean zero. Together with Lemma \ref{le:se}, this gives
     \[
     \E_{u^{[1]},\xi^{[1]}} [\|\varepsilon^0\|_p^2]
     =
     \E_{u^{[1]},\xi^{[1]}}[\|\frac{1}{n_0}\sum_{j=1}^{n_0}Z_j^0\|_p^2]
     \leq \frac{p-1}{n_0^2}\sum_{j=1}^{n_0}\E_{u^{[1]},\xi^{[1]}}[\|Z_j^0\|_p^2]
     \leq \frac{2(p-1)\tilde \sigma^2}{n_0}.
     \]
      Set
     \[
     A_k:=\E_{u^{[k+1]},\xi^{[k+1]}}[\|\varepsilon^k\|_p^2],
     \qquad
     D_k:=\E_{u^{[k+1]},\xi^{[k+1]}}
     [\mathcal{L}(x^k,\lambda^k;\mu)
     -\mathcal{L}(x^{k+1},\lambda^{k+1};\mu)]
     +m\rho_k\bar{M}_c^2,
     \]
     then (\ref{edx2}) can be represented as
     \[
     \alpha A_k\leq A_k-A_{k+1}+4\alpha^2  (p -1)\frac{\tilde \sigma^2}{n} + \frac{18(p-1)S_p\nu^2L_f^2\E_u[\|u\|_2^6]}{n}+ \frac{96 (p -1) \eta S_p L_f^2}{n}D_k.
     \]
     Summing above inequality over \(k=0,\ldots,K-2\) gives
     \begin{align*}
     \alpha\sum_{k=0}^{K-2}A_k
     &\leq A_0-A_{K-1}
     +\frac{4\alpha^2(p-1)\tilde\sigma^2(K-1)}{n}\\
     &\quad
     +\frac{18(p-1)S_p\nu^2L_f^2\E_u[\|u\|_2^6](K-1)}{n}
     +\frac{96(p-1)\eta S_pL_f^2}{n}
     \sum_{k=0}^{K-2}D_k.
     \end{align*}
     Adding \(\alpha A_{K-1}\) to both sides, since \(0<\alpha\leq1\), \(\sum_{k=0}^{K-2}D_k=\mathcal{L}(x^0,\lambda^0;\mu)
    -\E_{u^{[K-1]},\xi^{[K-1]}}[\mathcal{L}(x^{K-1},\lambda^{K-1};\mu)]
    +m \sum_{k=0}^{K-2}\rho_k \bar{M}_c^2,\) it follows from (\ref{del}) that
    \begin{align*}
         \alpha\sum_{k=0}^{K-1} A_k 
       &\leq  \frac{2(p-1)\tilde \sigma^2}{n_0}+\frac{4\alpha^2 (p -1)\tilde \sigma^2K}{n} + \frac{96 (p -1) \eta S_p L_f^2 \Delta_0}{n} 
         +\frac{18(p-1)S_p\nu^2L_f^2\E_u[\|u\|_2^6]K}{n}.
    \end{align*}
    Then dividing the whole inequality by $\alpha K$ yields 
    \be \label{evar}
    \begin{aligned}
         \frac{1}{K} \sum_{k=0}^{K-1} \E_{u^{[k+1]},\xi^{[k+1]}}[\|\varepsilon^k\|_p^2] \leq \frac{2(p-1)\tilde \sigma^2}{n_0\alpha K}+\frac{4\alpha (p -1)\tilde \sigma^2}{n} + \frac{3\nu^2\E_u[\|u\|_2^6]}{16\eta^2 } +\frac{ \Delta_0}{ \eta K}.
    \end{aligned}
    \ee
    Finally, with \(\tilde \lambda^k := \lambda^k+\mu c(x^k)\), combining the above results leads to
    \begin{align}
        & \E_{R;u^{[K]},\xi^{[K]}}  [\|\mathcal{G}_X (x^R,\nabla f(x^R) + \nabla c(x^R)^\top\tilde \lambda^R,\eta)\|_q^2 ]\notag\\
        &= \frac{1}{K} \sum_{k=0}^{K-1} \E_{u^{[k]},\xi^{[k]}} [\|\mathcal{G}_X (x^k,\nabla f(x^k) + \nabla c(x^k)^\top(\lambda^k+\mu c(x^k)),\eta)\|_q^2 ]\notag\\
        &\leq \frac{3}{K} \sum_{k=0}^{K-1} \E_{u^{[k]},\xi^{[k]}} [\|\nabla f(x^k) - g_\nu(x^k)\|_p^2 ] + \frac{3}{K} \sum_{k=0}^{K-1} \E_{u^{[k+1]},\xi^{[k+1]}} [\|g_\nu(x^k) - s^k \|_p^2 ]\notag\\
        &\quad+\frac{3}{K} \sum_{k=0}^{K-1} \E_{u^{[k+1]},\xi^{[k+1]}} [\|\mathcal{G}_X (x^k,s^k + \nabla c(x^k)^\top(\lambda^k+\mu c(x^k)),\eta)\|_q^2 ]\notag\\
        &\leq \frac{3\nu^2 L_f^2\E_u[\|u\|_2^6]}{4}+\frac{3}{K} \sum_{k=0}^{K-1} \E_{u^{[k+1]},\xi^{[k+1]}}[\|\varepsilon^k\|_p^2] +  \frac{3}{ K} \sum_{k=0}^{K-1} \frac{1}{\eta^2} \E_{u^{[k+1]},\xi^{[k+1]}} \|x^{k+1}-x^k\|_q^2 \label{eq:5.7}\\
         &\leq(4L_f^2+3\eta^{-2})\nu^2 \E_u[\|u\|_2^6] +\frac{30(p-1)\tilde \sigma^2}{n_0 \alpha K}+\frac{60\alpha (p -1)\tilde \sigma^2}{n} +\frac{27 \Delta_0}{ \eta K}, \notag
    \end{align}
    where the first inequality follows from the 1-Lipschitz continuity of $\mathcal{G}_X(x,\cdot,\eta)$ (see Lemma 6.5 and Proposition 6.1 in \cite{lan2020first}), the second inequality is due to Lemma  \ref{le:vg} and definitions of $\varepsilon^k$ and $\mathcal{G}_X$, {and the last inequality comes from \eqref{evar} and \eqref{edx}}. 
Therefore, \eqref{bound} is derived. 
\end{proof}

Consider the four terms on the right-hand side (R.H.S.) of \eqref{bound}. The first term comes from the bias of the zeroth-order gradient estimator, which is affected by the sixth-order moment \(\E_u[\|u\|_2^6]\) and can be controlled by choosing $\nu$  sufficiently small. The second and third terms are related to the variance $\tilde \sigma^2$ on \(G_\nu(x;u,\xi)\), and are controlled by the parameter $\alpha$ used for variance reduction. The last term depends on the gap between the value of the augmented Lagrangian function at the initial point and its lower bound and  generally decreases as the number of iterations grows.

We next concentrate on the primal feasibility.  

\begin{theorem}[{\bf {Primal feasibility}}]\label{th:3eKKT}
    Under conditions of Theorem \ref{stationary} and Assumption \ref{ass:ensc}, 
    %suppose that \red{\( v(x) \) is \( L_v \)-smooth within the constraint set \(X\) w.r.t. the \( \ell_q \)-norm.} 
    it holds that
    \be\label{bound-fea}\begin{aligned}
\mathbb{E}_{R;u^{[K]},\xi^{[K]}} [  \|  c ( x^R  )  \|_p^2]  
        \leq \frac{\|c(x^0)\|_p^2}{K} + \frac{8(M_f^2+M_h^2+m^2\rho^2\bar{M}_c^2M_c^2)}{\beta^2 \mu^2 }
        +\frac{8(8L_v^2+3)\Delta_0}{\beta^2\mu^2\eta K}\\+\frac{4(2L_v^2+1)\nu^2\E_u[\|u\|_2^6]}{\beta^2\mu^2\eta^2} +\frac{32(2L_v^2+1)(p-1)}{\beta^2\mu^2}(\frac{1}{n_0 \alpha K} + \frac{2\alpha}{n})\tilde \sigma^2.
    \end{aligned}\ee

\end{theorem}

\begin{proof}
  
    By the definition of 
    $\tilde{\lambda}^k = \lambda^k+\mu c(x^k)$, 
     there exists $g_h^k \in \partial h(x^k)$ such that
    \bee\begin{aligned}
        0 &= \, {\rm dist}_p ( s^{k-1}+g_h^k+ \nabla c ( x^{k-1} )^\top \tilde{\lambda}^{k-1}+\frac{1}{\eta}(\nabla v(x^k)-\nabla v(x^{k-1})),-\mathcal{N}_X(x^k)  )\\
        &=\, {\rm dist}_p (s^{k-1}+\partial h(x^k)+ \nabla c ( x^{k-1} )^\top \tilde{\lambda}^{k-1}+\frac{1}{\eta}(\nabla v(x^k)-\nabla v(x^{k-1})),-\mathcal{N}_X(x^k)  ).
    \end{aligned}\eee
    Then it holds from Assumption \ref{ass:ensc} and $\|\cdot\|_p\le \|\cdot\|_2 $ that for any  $k \in \{1,\ldots,K\}$,
    % \be \label{bof}
    \begin{align}
          & \|    c ( x^k )    \|_p^2 
         \le \, \frac{{\rm dist}_p^2( \nabla c ( x^k )^\top c ( x^k )  ,  -\mathcal{N}_{X}(x^k)  ) }{\beta^2} 
         = \frac{{\rm dist}_p^2 (\mu \nabla c ( x^k )^\top c ( x^k ) , -\mathcal{N}_X(x^k)  )}{\beta^2 \mu^2 }\notag \\
        & \leq  \frac{1}{\beta^2 \mu^2 }\|s^{k-1}+g_h^k+ \nabla c ( x^{k-1} )^\top \tilde{\lambda}^{k-1}+\frac{1}{\eta}(\nabla v(x^k)-\nabla v(x^{k-1}))-\mu \nabla c ( x^k )^\top  c ( x^k ) \|_p^2  \notag \\
        & \leq   \frac{8}{\beta^2 \mu^2 }(\|\nabla f(x^{k-1})\|_p^2+\|g_\nu(x^{k-1})-\nabla f(x^{k-1})\|_p^2+\|s^{k-1}-g_\nu(x^{k-1})\|_p^2+\|g_h^k\|_p^2)\notag \\
        &\quad +\frac{8L_v ^2}{\beta^2\eta^2\mu^2}\|x^k-x^{k-1} \|_q^2+\frac{8}{\beta^2\mu^2}\|\nabla c(x^{k-1})^\top\lambda^{k-1}\|_p^2\notag \\
        &\quad +\frac{8}{\beta^2}(\|\nabla c(x^{k-1})^\top(c(x^{k-1})-c(x^k))\|_p^2+\|(\nabla c(x^{k-1})-\nabla c(x^k))^\top c(x^k)\|_p^2)\notag \\
        &\leq  \frac{8}{\beta^2 \mu^2 }(M_f^2+M_h^2+\|g_\nu(x^{k-1})-\nabla f(x^{k-1})\|_p^2+\|\varepsilon^{k-1}\|_p^2)\notag \\
        &\quad +\frac{8L_v ^2}{\beta^2\eta^2\mu^2}\|x^k-x^{k-1} \|_q^2+\frac{8m}{\beta^2 \mu^2}\sum_{i=1}^m\|\lambda_i^{k-1}\nabla c_i(x^{k-1})\|_p^2\notag \\
        &\quad +\frac{8m}{\beta^2}\sum_{i=1}^m(\|(c_i(x^{k-1})-c_i(x^k))\nabla c_i(x^{k-1})\|_p^2+\|c_i(x^k)(\nabla c_i(x^{k-1})-\nabla c_i(x^k))\|_p^2)\notag \\
        &\leq \frac{8}{\beta^2 \mu^2 }(M_f^2+M_h^2+\|g_\nu(x^{k-1})-\nabla f(x^{k-1})\|_p^2+\|\varepsilon^{k-1}\|_p^2+m^2\rho^2\bar{M}_c^2M_c^2)\notag \\
        &\quad +\frac{8}{\beta^2}(\frac{L_v ^2}{\eta^2\mu^2}+m^2({M_c^4}+L_c^2\bar{M}_c^2))\|x^k-x^{k-1} \|_q^2, \notag 
    \end{align}
    where the third inequality follows from Jensen's inequality, the fourth one comes from $\|\nabla f(x)\|_2^2 \leq M_f^2$, $\|g_h^k\|_2^2 \leq M_h^2$, the last one is due to Cauchy-Schwarz inequality, Assumption \ref{ass:lip} and Lemma \ref{new:z2}. 
    Taking expectation w.r.t. $u^{[K]},~\xi^{[K]}$ and \(R\) on both sides of the above relation and using \eqref{evar} and \eqref{edx} yield
    \begin{align}
        &\E_{R;u^{[K]},\xi^{[K]}} [  \|    c ( x^R  )    \|_p^2 ]  
        =\frac{1}{K} ( \sum_{k=1}^{K-1} \E_{u^{[k]},\xi^{[k]}} [\ \|    c ( x^k  )    \|_p^2 ]  +\|c(x^0)\|_p^2 )\notag\\
        & \leq \frac{8}{\beta^2 \mu^2 }(M_f^2+M_h^2+\frac{\nu^2L_f^2\E_u[\|u\|_2^6]}{4}+m^2\rho^2\bar{M}_c^2M_c^2)+\frac{8}{\beta^2 \mu^2K} \sum_{k=1}^K \E_{u^{[k]},\xi^{[k]}}[\|\varepsilon^{k-1}\|_p^2]\notag\\
        &\quad +\frac{2(4L_v ^2+1)}{\eta^2\beta^2\mu^2K}\sum_{k=1}^K\E_{u^{[k]},\xi^{[k]}}[\|x^k - x^{k-1} \|_q^2] +\frac{\|c(x^0)\|_p^2}{K}\notag\\
        & \leq  \frac{\|c(x^0)\|_p^2}{K} + \frac{8}{\beta^2 \mu^2 }(M_f^2+M_h^2+\frac{\nu^2L_f^2\E_u[\|u\|_2^6]}{4}+m^2\rho^2\bar{M}_c^2M_c^2)\notag\\
        &\quad + \frac{2(4L_v^2+1)}{\beta^2\mu^2}(\frac{4\Delta_0}{\eta K}+\nu^2L_f^2\E_u[\|u\|_2^6])+\frac{16(2L_v^2+1)}{\beta^2\mu^2K}\sum_{k=1}^K \E_{u^{[k]},\xi^{[k]}}[\|\varepsilon^{k-1}\|_p^2]\label{eq:5.9}\\
        &\leq \frac{\|c(x^0)\|_p^2}{K} + \frac{8}{\beta^2 \mu^2 }(M_f^2+M_h^2+m^2\rho^2\bar{M}_c^2M_c^2)+\frac{8(8L_v^2+3)\Delta_0}{\beta^2\mu^2\eta K}\notag\\
        &\quad+\frac{4(2L_v^2+1)\nu^2\E_u[\|u\|_2^6]}{\beta^2\mu^2\eta^2} +\frac{32(2L_v^2+1)(p-1)}{\beta^2\mu^2}(\frac{1}{n_0 \alpha K} + \frac{2\alpha}{n})\tilde \sigma^2,\notag
    \end{align}
    where the first inequality uses \(\eta^2m^2({M_c^4}+L_c^2\bar{M}_c^2) \leq \frac{\eta^2L_\mu^2}{\mu^2} \leq \frac{1}{4\mu^2}\), the second inequality is due to \eqref{edx}, and the last inequality comes from \eqref{evar}.
\end{proof}

\subsection{Oracle complexity} \label{sec3-2} 
The upper bounds on the KKT measure depend on the parameter \(p\). As $p$ approaches infinity, the estimate will grow unbounded, rendering the estimate ineffective. To address this issue, we consider two regimes: \(p \in [2, 2 \ln d]\) and  \(p > 2 \ln d\). 
In the following discussion, let $\rho>0$ be a fixed constant and set
\(
    \rho_k=\frac{\rho}{K},  k=0,\ldots,K-1. 
\) 

We first establish the oracle complexity of Algorithm \ref{algo} for returning an $\epsilon$-KKT point of \eqref{p} when $p\in {[2,{2\ln d}]}$. In this case, we set
 the parameters $n_0,n,\nu,\mu,\alpha,\eta,K$ such that
% \bee
% \begin{gathered}
%     n_0 = \lceil\frac{5(p-1)\tilde \sigma^2}{4\epsilon^2}\rceil,n = \lceil(p-1)S_p\rceil,\, \nu^2 = \min\{\frac{\epsilon^2}{4(4L_f^2+3\eta^{-2})\E_u[\|u\|_2^6]},\frac{\beta^2\mu^2\eta^2\epsilon^2}{20(2L_v^2+1)\E_u[\|u\|_2^6]}\},\\
%     \mu^2 = \max\{\frac{10(8L_v^2+3)}{27\beta^2},\frac{8(2L_v^2+1)}{3\beta^2},\frac{40(M_f^2+M_h^2+m^2\rho^2\bar{M}_c^2M_c^2)}{\beta^2\epsilon^2}\},\\
%     \alpha = 96L_f^2\eta^2,\, \eta^2 = \min \{\frac{1}{4L_\mu^2},\frac{1}{96L_f^2},\frac{S_p\epsilon^2}{23040L_f^2\tilde \sigma^2}\}, K = \lceil\max \{\frac{108\Delta_0}{\eta\epsilon^2}, \frac{5\|c(x^0)\|_2^2}{\epsilon^2},\frac{1}{L_f^2\eta^2}\}\rceil.
% \end{gathered}
% \eee
% For an explicit and sequentially well-defined parameter choice, the apparent circular dependence among $\bar{\sigma}^2$, $\nu$, and $\eta$ in \eqref{settings} can be avoided by defining \(\bar{\sigma}^2 = \frac{\epsilon^2}{16} +4S_p(\sigma^2+M_f^2) +2M_f^2. \) Accordingly, the parameters can be chosen sequentially as 
\be\label{settings}
\begin{gathered}
\bar{\sigma}^2 = \frac{\epsilon^2}{16} +4S_p(\sigma^2+M_f^2) +2M_f^2, \quad \eta^2 = \min\{ \frac{1}{4L_\mu^2}, \frac{1}{96L_f^2}, \frac{S_p\epsilon^2} {23040L_f^2\bar{\sigma}^2} \},\\
\quad  \mu^2 = \max\{ \frac{10(8L_v^2+3)}{27\beta^2}, \frac{8(2L_v^2+1)}{3\beta^2}, \frac{40(M_f^2+M_h^2+m^2\rho^2\bar{M}_c^2M_c^2)} {\beta^2\epsilon^2} \},  \quad  n_0 = \lceil \frac{5(p-1)\bar{\sigma}^2}{4\epsilon^2} \rceil, \\ 
 \alpha=96L_f^2\eta^2, \quad   \nu^2 = \min\{ \frac{\epsilon^2} {4(4L_f^2+3\eta^{-2})\mathbb{E}_u[\|u\|_2^6]}, \quad  \frac{\beta^2\mu^2\eta^2\epsilon^2} {20(2L_v^2+1)\mathbb{E}_u[\|u\|_2^6]} \}, \\
 n=\lceil(p-1)S_p\rceil, \quad  K = \lceil \max\{ \frac{108\Delta_0}{\eta\epsilon^2}, \ \frac{5\|c(x^0)\|_2^2}{\epsilon^2}, \frac{1}{L_f^2\eta^2} \} \rceil.
 % \quad \rho_k=\frac{\rho}{K},  \ k=0,\ldots,K-1.
 \end{gathered}\ee
 
\begin{theorem}[{\bf {Oracle} complexity:  $p\in{[2,{2\ln d}]}$}]\label{thm:p}
Under Assumptions \ref{ass:lip}--\ref{ass:exs},  Assumption \ref{ass:v}--\ref{ass:ensc} and parameter setting \eqref{settings}, where all \(p\)-dependent assumptions are specified with \(p \in [2, 2\ln d]\). If \(\|c(x^0)\|_2^2\le 1/\mu\), then   Algorithm \ref{algo} returns an \(\epsilon\)-KKT point, with the oracle complexity of order  \(\mathcal{O}(pS_p\epsilon^{-3})\).
\end{theorem}
\begin{proof}
It follows from \eqref{settings} that \(\tilde{\sigma}^2= \nu^2 L_f^2 \E_u[\|u\|_2^6]+4 S_p (\sigma^2+M_f^2)+ 2M_f^2\leq \bar{\sigma}^2\).
We first verify the $\epsilon$-KKT conditions. 
By \eqref{bound} and the parameter setting \eqref{settings}, we obtain
\[
\begin{aligned}
&\E_{R;u^{[K]},\xi^{[K]}}
[
\|
\mathcal{G}_X
(
x^R,\nabla f(x^R)
+\nabla c(x^R)^\top\tilde{\lambda}^R,\eta
)
\|_q^2
]                                                     \\
&\leq
(4L_f^2+3\eta^{-2})\nu^2\E_u[\|u\|_2^6]
+\frac{30(p-1)\tilde{\sigma}^2}{n_0\alpha K}
+\frac{60\alpha(p-1)\tilde{\sigma}^2}{n}
+\frac{27\Delta_0}{\eta K}\\
&\leq \frac{\epsilon^2}{4}+\frac{\epsilon^2}{4}+\frac{\epsilon^2}{4}+\frac{\epsilon^2}{4}
= \epsilon^2
\end{aligned}
\]
We next verify primal feasibility.  By \eqref{bound-fea} and the parameter setting \eqref{settings}, we obtain
\[
\begin{aligned}
\mathbb{E}_{R;u^{[K]},\xi^{[K]}} [  \|  c ( x^R  )  \|_p^2]  
        &\leq \frac{\|c(x^0)\|_p^2}{K} + \frac{8(M_f^2+M_h^2+m^2\rho^2\bar{M}_c^2M_c^2)}{\beta^2 \mu^2 }
        +\frac{8(8L_v^2+3)\Delta_0}{\beta^2\mu^2\eta K}\\&+\frac{4(2L_v^2+1)\nu^2\E_u[\|u\|_2^6]}{\beta^2\mu^2\eta^2} +\frac{32(2L_v^2+1)(p-1)}{\beta^2\mu^2}(\frac{1}{n_0 \alpha K} + \frac{2\alpha}{n})\tilde \sigma^2\\
        &\leq \frac{\epsilon^2}{5}+\frac{\epsilon^2}{5}+\frac{\epsilon^2}{5}+\frac{\epsilon^2}{5}+\frac{\epsilon^2}{5}=\epsilon^2.
\end{aligned}
\]
Thus, $(x^R,\tilde{\lambda}^R)$ is an $\epsilon$-KKT point.
The number of stochastic function evaluations is bounded by
\begin{align} 
&  2n_0 + 4nK 
\leq \lceil\frac{5(p-1)\bar \sigma^2}{2\epsilon^2}\rceil + 4n\max \lceil\{\frac{108\Delta_0}{\eta\epsilon^2}, \frac{5\|c(x^0)\|_2^2}{\epsilon^2},\frac{1}{L_f^2\eta^2}\}\rceil \notag \\
&\leq \lceil\frac{5(p-1)\bar \sigma^2}{2\epsilon^2}\rceil + \lceil4(p-1)S_p\max \{\frac{216L_\mu\Delta_0}{\epsilon^2}, \frac{16394L_f\Delta_0\bar \sigma}{\sqrt{S_p}\epsilon^3}, \frac{5\|c(x^0)\|_2^2}{\epsilon^2},\frac{4L_\mu^2}{L_f^2}, 96, \frac{23040\bar \sigma^2}{S_p \epsilon^2}\}\rceil\notag \\
&= \mathcal{O}( \frac{pS_p\sigma^2}{\epsilon^2} + \frac{pS_pL_f\Delta_0}{\epsilon^2} + \frac{pS_p\tilde L_c\Delta_0}{\beta\epsilon^3} + \frac{pS_pL_f\Delta_0\sigma}{\epsilon^3}+\frac{pS_p\|c(x^0)\|_2^2}{\epsilon^2} + \frac{pS_p\tilde L_c^2}{\beta^2\epsilon^2}),\label{eq:complexity}
\end{align}
where the second inequality follows from the parameter choices of $\eta^2$, and the last line comes from $L_\mu = \mathcal O(L_f+\frac{\tilde L_c}{\beta \epsilon})$ with $\tilde L_c = M_c^2 + L_c$ and $\bar \sigma^2 = \mathcal O(S_p)$. We thus obtain the conclusion from  \(\Delta_0 = \mathcal L(x^0,\lambda^0;\mu) - \Phi_{\rm inf} + 2m\rho \bar{M}_c^2 = \mathcal O (\mu \|c(x^0)\|_2^2+1) = \mathcal O(1)\).
\end{proof}

\begin{remark}
The order \(\mathcal O(\epsilon^{-3})\) results from combining the
\(\mathcal O(\epsilon^{-1})\) penalty parameter with the stochastic
variance-reduction schedule.  The factor \(S_p\) captures the dimensional
effect of the zeroth-order estimator, while the constraint regularity condition plays a key role in translating approximate stationarity into primal feasibility.
\end{remark}

We next explore the oracle complexity when \(p > 2\ln d\). 
In this regime, the mirror descent geometry remains the original \(\ell_q\) geometry, where \(\frac{1}{p}+\frac{1}{q}=1\), while the stochastic estimator errors are controlled in the smaller norm exponent \(2\ln d\). Define
\[\bar p=2\ln d,\quad \tilde \sigma_{\bar p}^{\,2} = 
\nu^2 L_f^2 \E_u[\|u\|_2^6]+4 S_{\bar p} (\sigma^2+M_f^2)+ 2M_f^2,\quad 
\bar\sigma_{\bar p}^{\,2}
=\frac{\epsilon^2}{16}
+4S_{\bar p}(\sigma^2+M_f^2)+2M_f^2.
\]  %We next provide the detailed proof.

\begin{theorem}[\bf {Oracle} complexity:  {\(p>{2\ln d}\)}]\label{thm:pinf}
Under Assumptions \ref{ass:lip}--\ref{ass:ms-p} and Assumption \ref{ass:v}--\ref{ass:ensc} where all \(p\)-dependent assumptions are specified with  \(p > {2\ln d}\). Moreover, suppose Assumption \ref{ass:exs} holds at \(\bar p\) with constant \(S_{\bar p}\) and  parameters following  \eqref{settings} with \(p\), \(S_p\), \(\bar\sigma^2\) replaced by
\(\bar p\), \(S_{\bar p}\), 
\(\bar\sigma_{\bar p}^{\,2}\). If   \(\|c(x^0)\|_2^2\le1/\mu\), 
% in particular,
% \(n_0=\lceil\frac{5(\bar p-1)\bar\sigma_{\bar p}^{\,2}}{4\epsilon^2}\rceil,
% \ 
% n=\lceil(\bar p-1)S_{\bar p}\rceil,\ 
% \eta^2 = \min \{\frac{1}{4L_\mu^2},\frac{1}{96L_f^2},
% \frac{S_{\bar p}\epsilon^2}
% {23040L_f^2\bar\sigma_{\bar p}^{\,2}}\},\)
then Algorithm \ref{algo} returns an
\(\epsilon\)-KKT point with oracle complexity in order 
\(\mathcal{O}(S_{2\ln d}\ln d\cdot\epsilon^{-3})\).
\end{theorem}
\begin{proof}
The proof idea is essentially the same as that of Theorem \ref{thm:p}. It only requires replacing, respectively, the stationary upper bound under the \(\ell_q\)-norm in Theorem \ref{stationary} and the primal feasibility upper bound under the \(\ell_p\)-norm in Theorem \ref{th:3eKKT} by their corresponding \(\ell_{\bar p}\)-type quantities that serve as upper bounds.  

First, we give an upper bound on the stationarity. Since \(d\ge3\), we have \(\bar p=2\ln d\ge2\).  Moreover,
due to \(p\ge\bar p\) and  the monotonicity of finite-dimensional
\(\ell_r\)-norm as a function of \(r\), it gives \( \|z\|_p\le \|z\|_{\bar p}\)
for any \(    z\in\mathbb R^d.\)
% \begin{equation}\label{eq:pinf-norm}
%     \|z\|_p\le \|z\|_{\bar p},
%     \qquad z\in\mathbb R^d.
% \end{equation} 
Since
\(K\ge \frac{1}{L_f^2\eta^2}\)
and
\(\eta^2\le\frac{1}{96L_f^2},
\)
we have \(K\ge2\). Set
\[
\bar A_k
=\E_{u^{[k+1]},\xi^{[k+1]}}
  [\|\varepsilon^k\|_{\bar p}^2], \quad  k=0,\ldots,K-1.
\] 
As Assumption \ref{ass:exs} holds at \(\bar p\), the arguments
of Lemmas \ref{le:vg}, \ref{le:lip}, and \ref{le:se} apply with
\(p\) replaced by \(\bar p\).
In particular, repeating the proof of Lemma~\ref{le:eps} and using
the \((\bar p-1)\)-smoothness of
\(\frac12\|\cdot\|_{\bar p}^2\), we obtain
\begin{align*}
\bar A_{k+1}
\leq
(1-\alpha)^2\bar A_k
&+\frac{4\alpha^2(\bar p-1)}{n}
    \tilde\sigma_{\bar p}^{\,2} +
\frac{12(\bar p-1)(1-\alpha)^2}{n}
    \nu^2L_f^2\E_u[\|u\|_2^6]\\
&+
\frac{24(\bar p-1)(1-\alpha)^2S_{\bar p}L_f^2}{n}
    \E[\|x^{k+1}-x^k\|_2^2].   
\end{align*}
Furthermore, since
\(\eta^2\le\frac{1}{4L_\mu^2}\), we have
\(\eta L_\mu\le\frac{1}{2}\). Hence, as in the proof of
Theorem~\ref{stationary},
\[
\frac1K\sum_{k=0}^{K-1}
\E[\|x^{k+1}-x^k\|_q^2]
\leq{}
\frac{4\eta\Delta_0}{K}
+\eta^2\nu^2L_f^2\E_u[\|u\|_2^6]  +
\frac{4\eta^2}{K}
\sum_{k=0}^{K-1}\bar A_k,
\]
where %we have also used \eqref{eq:pinf-norm} to bound
\(\|\varepsilon^k\|_p^2\le \|\varepsilon^k\|_{\bar p}^2\) due to $p\ge \bar p$.
Since
\( n\ge(\bar p-1)S_{\bar p},
    \alpha=96L_f^2\eta^2,
\)
we have
\(\frac{96(\bar p-1)S_{\bar p}\eta^2L_f^2}{\alpha}
    \le n.\)
Moreover,
\(\eta^2\le\frac{1}{96L_f^2}\) implies \(0<\alpha\le1\).
Repeating the summation argument in the proof of
Theorem~\ref{stationary}  yields
\begin{equation}\label{eq:pinf-evar}
\frac1K\sum_{k=0}^{K-1}\bar A_k
\leq
\frac{2(\bar p-1)\tilde\sigma_{\bar p}^{\,2}}
     {n_0\alpha K}
+\frac{4\alpha(\bar p-1)\tilde\sigma_{\bar p}^{\,2}}{n}  
+
\frac{3\nu^2\E_u[\|u\|_2^6]}{16\eta^2}
+\frac{\Delta_0}{\eta K}.
\end{equation}
Note that 
 \(\E[\|\varepsilon^k\|_p^2]
    \leq    \E[\|\varepsilon^k\|_{\bar p}^2]
    =\bar A_k.\)
And applying Lemma~\ref{le:vg} at exponent \(\bar p\) gives
\[
    \|g_\nu(x^k)-\nabla f(x^k)\|_p^2
    \leq
    \|g_\nu(x^k)-\nabla f(x^k)\|_{\bar p}^2
    \leq
    \frac{\nu^2L_f^2}{4}\E_u[\|u\|_2^6].
\]
Consequently, following the same decomposition as in
\eqref{eq:5.7}, we obtain that with \(\tilde\lambda^R
    =
    \lambda^R+\mu c(x^R),\)
\begin{align}
&\E_{R;u^{[K]},\xi^{[K]}}
[\|
\mathcal G_X
(
x^R,
\nabla f(x^R)+
\nabla c(x^R)^\top\tilde\lambda^R,
\eta)\|_q^2]                                                \notag\\
&\leq
(4L_f^2+3\eta^{-2})
\nu^2\E_u[\|u\|_2^6]
+\frac{30(\bar p-1)\tilde\sigma_{\bar p}^{\,2}}
       {n_0\alpha K}         
+\frac{60\alpha(\bar p-1)\tilde\sigma_{\bar p}^{\,2}}{n}
+\frac{27\Delta_0}{\eta K}.
\label{eq:pinf-sta}
\end{align}
%where \(\tilde\lambda^R = \lambda^R+\mu c(x^R).\)
Then substituting the new parameter settings into \eqref{eq:pinf-sta} and applying  \(\tilde \sigma_{\bar p}^{\,2} \leq \bar\sigma_{\bar p}^{\,2}\) yield
\[\E
[
\|\mathcal G_X(
 x^R,\nabla f(x^R)+\nabla c(x^R)^\top\tilde\lambda^R,\eta
)\|_q^2
]\le
\frac{\epsilon^2}{4}
+\frac{\epsilon^2}{4L_f^2\eta^2K}
+\frac{\epsilon^2}{4}
+\frac{27\Delta_0}{\eta K}
\le\epsilon^2,
\]
where the final inequality follows from
\( K\ge
    \max \{
        \frac{108\Delta_0}{\eta\epsilon^2},
        \frac{1}{L_f^2\eta^2}
     \}.\)

We next verify primal feasibility. In the proof of Theorem~\ref{th:3eKKT}, the terms
\(\|g_\nu(x^k)-\nabla f(x^k)\|_p^2\) and \(\|\varepsilon^k\|_p^2\)
can both be bounded by their corresponding
\(\ell_{\bar p}\)-norm quantities due to $p
\ge \bar p.$ 
Therefore, using the same argument as in Theorem~\ref{th:3eKKT},  substituting the new parameter settings and applying \(\tilde \sigma_{\bar p}^{\,2} \leq \bar\sigma_{\bar p}^{\,2}\) derive
\[\begin{aligned}
\mathbb{E}_{R;u^{[K]},\xi^{[K]}} [\|c(x^R)\|_p^2]  
        &\leq \frac{\|c(x^0)\|_p^2}{K} + \frac{8(M_f^2+M_h^2+m^2\rho^2\bar{M}_c^2M_c^2)}{\beta^2 \mu^2 }
    +\frac{8(8L_v^2+3)\Delta_0}{\beta^2\mu^2\eta K}+\frac{4(2L_v^2+1)\nu^2\E_u[\|u\|_2^6]}{\beta^2\mu^2\eta^2} \\&\qquad+\frac{32(2L_v^2+1)(\bar{p}-1)}{\beta^2\mu^2}(\frac{1}{n_0 \alpha K} + \frac{2\alpha}{n})\tilde \sigma_{\bar{p}}^2\leq \frac{\epsilon^2}{5}+\frac{\epsilon^2}{5}+\frac{\epsilon^2}{5}+\frac{\epsilon^2}{5}+\frac{\epsilon^2}{5}=\epsilon^2.
    \end{aligned}
\]
Hence, \(x^R\) is an
\(\epsilon\)-KKT point following Definition \ref{de:kkt}.

Under the given parameter settings,
\(\mu=\mathcal O(\epsilon^{-1}),
L_\mu=\mathcal O(\epsilon^{-1}),
\bar\sigma_{\bar p}^{\,2}
=\mathcal O(S_{\bar p}),
\eta^{-1}=\mathcal O(\epsilon^{-1}).\)
Furthermore, since \(\lambda^0=0\) and
\(\|c(x^0)\|_2^2\le1/\mu\), we have
\(\Delta_0
=\mathcal L(x^0,\lambda^0;\mu)-\Phi_{\rm inf}+2m\rho \bar{M}_c^2
=\mathcal O(1).\)
Hence,
\(K
=\lceil
\max\{
\frac{108\Delta_0}{\eta\epsilon^2},
\frac{1}{L_f^2\eta^2}\}
\rceil
=\mathcal O(\epsilon^{-3}).\)
Finally, the number of stochastic function evaluations is
bounded by
\[
\begin{aligned}
2n_0+4nK
&\le\lceil\frac{5(\bar p-1)\bar\sigma_{\bar p}^{\,2}}{2\epsilon^2}\rceil+\lceil4(\bar p-1)S_{\bar p}\max\{\frac{108\Delta_0}{\eta\epsilon^2},\frac{1}{L_f^2\eta^2}\}\rceil                 =\mathcal O(\bar pS_{\bar p}\epsilon^{-2}+\bar pS_{\bar p}\epsilon^{-3})   ,
\end{aligned}
\]
which is in order $\mathcal O(S_{2\ln d}\ln d\cdot\epsilon^{-3})$, thus deriving the conclusion. 
\end{proof}

\subsection[The mirror descent geometry and the value of Sp]{The mirror descent geometry and the value of $S_p$}
\label{ssec:cpx} 
In this subsection, we give discussions on  the mirror descent geometry and the values of \(S_p\).

{\bf Mirror descent geometry.}
In optimization algorithm's design, the choice of the mirror descent metric plays a pivotal role in determining the algorithm's effectiveness, particularly in how well it captures the underlying geometric structure of the problem. 
We explore the setup of the mirror descent framework, with an emphasis on selecting an appropriate distance-generating function based on the specific characteristics of the optimization problem under consideration.
The most common choice reduces to standard Euclidean geometry and arises from the distance-generating function  \( v(x) = \frac{1}{2} \|x\|_2^2 \). 
The distance-generating function is \(1\)-strongly convex and $1$-smooth over $X$ w.r.t. \(\ell_2\)-norm. 
However, for problems with special structures, such as some non-Euclidean geometry, alternative {mirror descent metrics} may better capture the intrinsic geometry and enhance the algorithm's performance. 

In  our algorithm, we require the distance-generating function satisfying Assumption \ref{ass:v}.  For \(p=2\), we use the standard Euclidean distance-generating function \( v(x) = \frac{1}{2} \|x\|_2^2 \). For \(p>2\), equivalently \(q\in(1,2)\), we employ a smoothed distance-generating function:
\begin{equation}\label{eq:smoothed_v}
v_\delta(x) = \frac{1}{2(q-1)} ( \sum_{i=1}^d (x_i^2 + \delta^2)^{q/2} )^{2/q},
\end{equation}
where \(\delta > 0\) is fixed independently of \(d,q\), and the target accuracy \(\epsilon\).
The incorporation of \(\delta\) avoids the Hessian singularities associated with zero coordinates in the standard squared \(\ell_q\)-norm. It is proved in Lemma  \ref{lem:smooth-mirror-map} that $v_\delta$ is  \(1\)-strong convex and $L_v$-smooth over $X$ w.r.t.\ the \( \ell_q \)-norm, where \(L_v\) is a constant dependent on \(q,d,\delta,D_X\) and defined by 
\begin{equation}\label{eq:Lv}
    L_{v}= \frac{2-q}{q-1} +
    \frac{1}{q-1}
    (D_X^q+d\delta^q)^{(2-q)/q}\delta^{q-2}.
\end{equation}

{\bf Value of $S_p$.} 
As can be seen from Theorems \ref{thm:p} and \ref{thm:pinf}, the value of \(S_p\) influences the dimensional dependence of the oracle complexity, with the oracle complexity scaling linearly with \(S_p\) due to \(n = (p-1)S_p\) for \(p \in [2, 2\ln d]\) and with \(S_{2 \ln d}\) due to \(n = (2 \ln d - 1)S_{2 \ln d}\) for \(p > 2 \ln d\). In spirit, the value of $S_p$ is determined by the probability distribution $\mathbb P_u$ on the random direction $u.$
%We now analyze the impact of different probability distributions on \( S_p \) under \( \ell_p \) norm. 

For \( p = 2 \), it is the case studied in most related work and the dimensional dependence is the least favorable among the cases studied in this paper. Note that 
$
    \E [\|\langle g,u\rangle u\|_2^2] = \E [g^\top uu^\top g u^\top u].
$ 
If  $u$ follows the {\it Rademacher distribution}  (corresponding to  Rademacher smoothing) or a {\it uniform distribution} on a sphere with the radius $\sqrt{d}$, we have $u^\top u = d$ and $\E [uu^\top] = {\rm I}_d$. Hence, $\E [\|\langle g,u\rangle u\|_2^2] = d \|g\|_2^2$, so that $S_2 = d$. If $u$ follows the {\it standard normal distribution}, i.e., $u_i, i=1,\ldots,d$, are independent standard Gaussian random variables, then $\E[u_i^2]=1$, $\E[u_i^4]=3$, and $\E[u_i u_j]=0$ for any $i \neq j$. Therefore,
\(
 \E[uu^\top uu^\top]=\E[\|u\|_2^2uu^\top]=(d+2)\mathrm I_d,
\) 
which gives $\E [\|\langle g,u\rangle u\|_2^2] = (d+2) \|g\|_2^2$. Thus, $S_2 = d+2$ in the Gaussian case. Consequently, for all three commonly used sampling distributions above, Theorem~\ref{thm:p} yields the oracle complexity  $\mathcal{O}(2S_2 \epsilon^{-3}) = \mathcal{O}(d \epsilon^{-3})$. Moreover, this linear dependence on d cannot be improved when \(p=2\), regardless of the choice of sampling distribution. Indeed, Theorem 2.2 of \cite{ma2025revisiting} shows that
\[
\mathbb{E}[\|\langle g, u \rangle u\|_2^2] \geq d\|g\|_2^2, \quad \forall g\in\R^d,
\]
This lower bound motivates the use of a non-Euclidean mirror descent geometry with \(p>2\), for which the dimensional dependence of the zeroth-order estimator can be reduced.  

For \( p > 2 \), the constant  \( S_p \) in Assumption \ref{ass:exs}is determined through \( \mathbb{E}[\|\langle g, u \rangle u\|_p^2] = \mathbb{E}[|\langle g, u \rangle|^2 \|u\|_p^2] \). For the standard normal distribution or the uniform distribution on the sphere of radius \( \sqrt{d} \), obtaining an explicit expression for \(S_p\)  is less straightforward. We thus focus on Rademacher smoothing. In this case, we have 
\begin{equation}\label{eq:Sp-rad}
\mathbb{E}[|\langle g, u \rangle|^2 \|u\|_p^2] = d^{2/p} \mathbb{E}[g^\top u u^\top g] = d^{2/p} \|g\|_2^2,
\qquad\text{and hence}\qquad S_p=d^{2/p}.
\end{equation}
Moreover, 
\(\E[\|u\|_2^6]=d^3<\infty.\)
Together with $\E[uu^\top]={\rm I}_d$, it verifies that under  Rademacher
smoothing Assumption~\ref{ass:exs} holds with
$S_p=d^{2/p}$.

\section{Oracle complexity under Rademacher smoothing}\label{sec:compl}

We will present the oracle complexity of Algorithm \ref{algo} under Rademacher smoothing, which gives \eqref{eq:Sp-rad}.

\subsection{Oracle complexity in mirror descent geometry}
We first give the main complexity  result in mirror descent geometry.
%For the smoothed distance-generating function defined in \eqref{eq:smoothed_v}, with $\frac{1}{p}+\frac{1}{q}=1$,its smooth Lipschitz constant (see Lemma \ref{lem:smooth-mirror-map}) is where $R$ is the $\ell_q$-radius of $X$ in Assumption \ref{ass:lip}.  
\begin{theorem}%[Oracle complexity under Rademacher smoothing]
\label{thm:rad_qnorm}
Let Algorithm \ref{algo} be equipped with the mirror descent metric
generated by \(v_\delta\) in \eqref{eq:smoothed_v}. Suppose that \(u\) follows the Rademacher
distribution and Assumptions \ref{ass:lip}--\ref{ass:ms-p} and Assumtion \ref{ass:ensc} hold. Further suppose that the bound $D_X$ in Assumption~\ref{ass:lip} and the regularity constant $\beta$ in Assumption~\ref{ass:ensc} are  independent of $d$ and $p$. 
Under the parameter setting following Theorem \ref{thm:p} if
\(p\in[2,2\ln d]\) and following  Theorem \ref{thm:pinf} if
\(p>2\ln d\),  if \(\|c(x^0)\|_2^2\leq 1/\mu,\)
then Algorithm \ref{algo} returns an \(\epsilon\)-KKT point, with the oracle complexity in order
\be\label{R-sm}
\begin{cases}
\mathcal O(p d^{2/p}\epsilon^{-3}
        +
        p^2 d\,\epsilon^{-2}
        +
        p^3 d^{2-2/p}), & \quad p\in[2,2\ln d],\\
\mathcal O(\ln d\cdot\epsilon^{-3}
        +
        p d^{1-2/p}\ln d\cdot\epsilon^{-2} +
        p^2 d^{2-4/p}\ln d),& \quad p>2\ln d.
\end{cases}
\ee

% \[\mathcal O(p d^{2/p}\epsilon^{-3}
%         +
%         p^2 d\,\epsilon^{-2}
%         +
%         p^3 d^{2-2/p}).\]
% When \(p>2\ln d\), the oracle complexity is \[\mathcal O(\ln d\cdot\epsilon^{-3}
%         +
%         p d^{1-2/p}\ln d\cdot\epsilon^{-2}
%         \\+
%         p^2 d^{2-4/p}\ln d
%     ).\]
\end{theorem}

\begin{proof}
By Lemma~\ref{lem:smooth-mirror-map}, the required Assumption  \ref{ass:v} on Theorems~\ref{thm:p} and \ref{thm:pinf}  is satisfied. Therefore, under the corresponding
parameter settings, Algorithm~\ref{algo} returns an
\(\epsilon\)-KKT point. Note from \eqref{eq:Lv} that \(L_v\) is a constant dependent on \(q,d\). So  the dependence must be retained in the oracle bound because $L_v$ enters the choices of $\mu$, $L_\mu$, $\eta$, and hence $K$.

We first consider the case \(p\in[2,2\ln d]\).
Under Rademacher smoothing, \eqref{eq:Sp-rad} gives
\(S_p=d^{2/p}\).
By the choice of the zeroth-order smoothing parameter \(\nu\) in
\eqref{settings}, we have
\(\nu^2L_f^2\E_u[\|u\|_2^6]\leq\epsilon^2/16\).
Since \(\beta\) is independent of \(d\) and \(p\), from the choice of \(\mu\) in \eqref{settings}, we have
\(\mu^2=\mathcal O(\epsilon^{-2}+L_v^2)\).
By Lemma \ref{le:Lu}, this further gives
\(L_\mu=\mathcal O(\epsilon^{-1}+L_v)\).
Together with
\(\bar\sigma^2=\mathcal O(S_p)\), the definition of
\(\eta\) in \eqref{settings} yields
\(\eta^{-1}=\mathcal O(\epsilon^{-1}+L_v)\).
Furthermore, 
\(\Delta_0
=\mathcal L(x^0,\lambda^0;\mu)-\Phi_{\rm inf}+2m\rho \bar{M}_c^2
=\mathcal O(\mu\|c(x^0)\|_2^2+1)
=\mathcal O(1)\).
Therefore, it follows from the definition of \(K\) in \eqref{settings} that 
\be\label{K-or}
K=\mathcal O(
    \frac{1}{\eta\epsilon^2}
    +\frac{\|c(x^0)\|_2^2}{\epsilon^2}
    +\frac{1}{\eta^2}
)=\mathcal O(\epsilon^{-3}+L_v\epsilon^{-2}+L_v^2).\ee
 Since \(q=p/(p-1)\), $D_X$ is independent of $d$ and $p$ and $\delta>0$ is fixed, we rewrite \(L_v\) in
\eqref{eq:Lv}  as
\(
L_v
=
(p-2)
+
(p-1)
[
    (\frac{D_X}{\delta})^{p/(p-1)}
    +d
]^{1-2/p}=\mathcal O(p d^{1-2/p}),
\) 
which leads to 
\be\label{K-or-1}
K
=
\mathcal O(
    \epsilon^{-3}
    +p d^{1-2/p}\epsilon^{-2}
    +p^2d^{2-4/p}
).
\ee
On the other hand, since
\(\bar\sigma^2=\mathcal O(S_p)\) and
\(S_p=d^{2/p}\), the batch sizes satisfy
\(n_0=\mathcal O(p d^{2/p}\epsilon^{-2})\) and
\(n=\mathcal O(p d^{2/p})\).
Therefore, the total number of stochastic function value evaluations
satisfies
\[
\begin{aligned}
2n_0+4nK
&=
\mathcal O(
    p d^{2/p}\epsilon^{-2}
)
+
\mathcal O(
    p d^{2/p}
    (
        \epsilon^{-3}
        +p d^{1-2/p}\epsilon^{-2}
        +p^2d^{2-4/p}
    ))                                          =
\mathcal O(
    p d^{2/p}\epsilon^{-3}
    +p^2d\,\epsilon^{-2}
    +p^3d^{2-2/p}).
\end{aligned}
\]

We next consider the case \(p>2\ln d\). As established in Theorem \ref{thm:pinf},  under Rademacher smoothing,
\(S_{\bar p}=d^{2/\bar p}=d^{1/\ln d}=e\). Then it follows from the setting required by Theorem \ref{thm:pinf} that
\(\bar\sigma_{\bar p}^{\,2}
=\mathcal O(S_{\bar p})=\mathcal O(1), n_0=\mathcal O(\ln d\cdot\epsilon^{-2})\) and
\(n=\mathcal O(\ln d)\).
The quantities
\(\mu\), \(L_\mu\), and \(\eta\) retain their dependence on  \(L_v\). By the preceding argument \eqref{K-or} is derived. 
% \(K=\mathcal O(\epsilon^{-3}+L_v\epsilon^{-2}+L_v^2)\).
Using again \(L_v=\mathcal O(p d^{1-2/p})\), we obtain \eqref{K-or-1}. 
% \[
% K
% =
% \mathcal O(
%     \epsilon^{-3}
%     +p d^{1-2/p}\epsilon^{-2}
%     +p^2d^{2-4/p}).
% \]
Finally,  the total number of stochastic function value evaluations is 
\[
\begin{aligned}
2n_0+4nK
&=
\mathcal O(
    \ln d\cdot\epsilon^{-2}
)
+
\mathcal O(
    \ln d
    (
        \epsilon^{-3}
        +p d^{1-2/p}\epsilon^{-2}
        +p^2d^{2-4/p}
    )
)                                                     \\
&=
\mathcal O(
    \ln d\cdot\epsilon^{-3}
    +p d^{1-2/p}\ln d\cdot\epsilon^{-2}
    +p^2d^{2-4/p}\ln d
).
\end{aligned}
\]
Hence, the conclusions are derived.
\end{proof}

The oracle complexity bounds in Theorem~\ref{thm:rad_qnorm}
reveal a trade-off in the choice of the mirror-geometry parameter
\(p\). 

For \(p\in[2,2\ln d]\), 
the first term \(p d^{2/p}\epsilon^{-3}\) originates from the dimension-dependent second moment
of the zeroth-order gradient estimator. Since
\(S_p=d^{2/p}\) under Rademacher smoothing, increasing \(p\)
reduces the dimensional factor \(d^{2/p}\) in this term.
The second and third terms arise from the smoothness constant of the
non-Euclidean mirror map. In particular, for the smoothed
distance-generating function \(v_\delta\),
\(L_v=\mathcal O(p d^{1-2/p})\), so the interaction between the
zeroth-order variance and the mirror smoothness produces the term
\(p^2d\,\epsilon^{-2}\), while the quadratic dependence on the
mirror smoothness produces \(p^3d^{2-2/p}\).
Hence, increasing \(p\) decreases the dimensional exponent in the
zeroth-order variance term but simultaneously increases the
geometry-dependent contribution.

For \(p>2\ln d\), the stochastic estimator is controlled using
\(\bar p=2\ln d\), so its batch-size contribution no longer improves
with further increase in \(p\). However, the dimensional exponents of the two geometry-dependent terms in the oracle complexity are \(1-2/p\) and \(2-4/p\),
respectively. The former is strictly smaller than \(1\) for every
\(p>2\), whereas the latter is smaller than \(1\) if and only if
\(p<4\). Hence, both exponents are strictly sublinear when
\(2\ln d<p<4.\)
Such an interval is nonempty only when \(d<e^2\). Therefore, this
regime may yield sublinear polynomial exponents in low dimensions,
but it does not provide a uniformly sublinear dimensional dependence
in the high-dimensional regime.
%\end{itemize}

An additional feature of the complexity bounds in \eqref{R-sm} is that the three
terms for both cases have different dependencies on the target accuracy
\(\epsilon\). The zeroth-order variance term scales
as \(\epsilon^{-3}\), whereas the two geometry-dependent terms scale
as \(\epsilon^{-2}\) and \(\epsilon^0\), respectively. In particular,
%even when the latter terms do not have sublinear polynomial dependence on \(d\), 
the variance term  dominates the oracle complexity when
%in the {\it high-accuracy regime}
\begin{equation}
\epsilon\le\min\{
   \frac{1}{p d^{1-2/p}},
    \frac{1}{p^{2/3}d^{\,2/3-4/(3p)}}
\}.
\tag{high-accuracy regime}
\end{equation} 
For \(p\in[2,2\ln d]\), this high-accuracy regime  
ensures that \(p d^{2/p}\epsilon^{-3}\) dominates both \(p^2d\,\epsilon^{-2}\) and
\(p^3d^{2-2/p}\), which implies the leading high-accuracy oracle complexity in order 
\(\mathcal O(p d^{2/p}\epsilon^{-3})\).
Similarly, when \(p>2\ln d\), the term \(\ln d\cdot\epsilon^{-3}\) dominates the two
geometry-dependent terms, yielding the leading order 
\(\mathcal O(
        \ln d\cdot\epsilon^{-3}).\)

\subsection{Oracle complexity in Euclidean geometry}
We now turn to the most common setting in the literature, where the Euclidean geometry is considered. 
We will prove that in this scenario the oracle complexity  exactly matches the oracle complexity bound, $\mathcal{O}(d \epsilon^{-3})$, of state-of-the-art {\it unconstrained}  stochastic zeroth-order optimization algorithms under the standard Euclidean setup~\cite{BG22,spider, huang2020accelerated, JWZL,nguyen2022stochastic}.

\begin{corollary}\label{thm:2-norm}
Let Algorithm \ref{algo} be equipped with the mirror descent metric generated by \(v(x) = \frac{1}{2}\|x\|_2^2\). Suppose that $u$ follows the Rademacher distribution,  Assumptions \ref{ass:lip}--\ref{ass:ms-p} and \eqref{CQ_2} hold. Then under   parameter setting \eqref{settings} with $p=2$, $S_2=d$, and $L_v=1$, 
if \(\|c(x^0)\|_2^2\le 1/\mu\), the oracle complexity of Algorithm \ref{algo} is of order  \(  \mathcal{O}(d\epsilon^{-3}) \) 
to find a point \(x\) satisfying \eqref{eps-KKT-2}. 
\end{corollary}

\begin{proof}
The result follows directly from Theorem \ref{thm:p} under the stated specialization of \eqref{settings}. %\blue{(the corollary  is Theorem 3 applied to \(p=2\))}
\end{proof}

Few stochastic zeroth-order methods in the literature directly address problem \eqref{p}. Wang et al.~\cite{wang2017penalty} propose a double-loop penalty method, but it incurs relatively high oracle complexity. The studies \cite{nguyen2022stochastic,SGG} consider settings in which only (stochastic) function values of the constraint functions are available and establish the corresponding complexity results. {Since \cite{SGG} relies on additional assumptions which are not imposed in our setting, such as the unbiasedness of the zeroth-order gradient estimators,} we do not compare against it here. Specializing the results of \cite{nguyen2022stochastic} to the case in which exact gradients of the constraint functions are available yields a complexity of \(\mathcal{O}(d \epsilon^{-4})\).
For stochastic first-order methods addressing problems with deterministic constraints, the state-of-the-art complexity, leveraging variance reduction techniques, is \(\mathcal{O}(\epsilon^{-3})\), as established in \cite{SWW,lu2026variance}, under the mean-squared smoothness assumption.
In our work, the oracle complexity order w.r.t. $\epsilon$ matches that in \cite{SWW,lu2026variance}, thereby improving on the complexity of existing constrained stochastic zeroth-order methods.

\subsection{Oracle complexity under a multi-stage scheme}
In the preceding analysis, all complexity results are based on the assumption of near-feasibility of the initial point. This requirement  stems from the need to ensure \(\Delta_0=\mathcal O(1)\) while simultaneously using a sufficiently large penalty parameter \(\mu\) to guarantee final feasibility. In general, however, such an initial point may not be readily available. A  strategy that gradually increases \(\mu\) is considered in \cite{alacaoglu2024complexity}, but with an extra factor $\log (\epsilon^{-1})$ in final  complexity bounds. Differently, we develop a multi-stage scheme that removes the feasibility of the starting point while preserving the leading high-accuracy oracle complexity bounds as the single-stage method. More specifically,  we use a more rapidly increasing penalty sequence so that only the first stage may have an \(L_v\)-dependent initial augmented-Lagrangian gap, while all subsequent stages admit a uniformly bounded initial gap. 
The corresponding method is summarized in Algorithm~\ref{algo:restart}. 

\begin{algorithm}[ht]
\caption{Multi-stage mirror descent linearized augmented Lagrangian method}
\label{algo:restart}
\begin{algorithmic}[1]
\REQUIRE \(\tilde x^0\in X\), \(0<\epsilon<1, \mu_{\rm base}, \gamma, \mu_*>0\)
\ENSURE $\tilde x^T$
\STATE{Set \(x_1^0=x^0\), \(\mu_1=\mu_{\rm base}\),
\(\mu_\star=\max\{\mu_{\rm base},\gamma/\epsilon\}\).}
%\STATE{Define \(\mu_{s+1}=\min\{2\mu_s^2,\mu_\star\}\) for \(s\ge1\), let \(T\) be the smallest index s.t. \(\mu_T=\mu_\star\).}
\FOR{\(s=1,\cdots\)}
\STATE{Set \(\epsilon_s=\gamma/\mu_s\) and \(\lambda_s^0=0\). Choose the stagewise parameters \(\nu_s\), \(\eta_s\), \(\alpha_s\), \(n_{0,s}\),\(n_s\) and \(K_s\).}
\STATE Run Algorithm~\ref{algo} from \(x_s^0\) with penalty parameter  \(\mu_s\), iteration number \(K_s\) and parameters \(\nu_s\), \(\eta_s\), \(\alpha_s\), \(n_{0,s}\), \(n_s\);   denote the output  by \(\tilde x^s\).
\IF{\(\mu_s \ge \mu_*\)}
\STATE{Terminate and return $\tilde x^T$ with $T:=s.$}
\ENDIF
%\IF{\(s<T\)}
\STATE{Set \(x_{s+1}^0=\tilde x^s\) and $\mu_{s+1}=\min\{2\mu_s^2,\mu_\star\}$.}
%\ENDIF
\ENDFOR
%\UNTIL{\(\mu_s=\mu_\star\)}\STATE{Set \(T=s\)}
%\STATE{\textbf{Return:} \(\tilde x^T\).}
\end{algorithmic} 
\end{algorithm}
%  At stage \(s\), Algorithm~\ref{algo} is initialized at \(x_s^0\) with \(\lambda_s^0=0\) and uses the penalty parameter \(\mu_s\). For \(0<\epsilon\le1\), set
% \begin{equation}\label{eq:restart-mu}
% \mu_\star=\max\{\mu_{\rm base},\frac{\gamma}{\epsilon}\},
% \qquad
% \mu_1=\mu_{\rm base},
% \qquad
% \mu_{s+1}=\min\{2\mu_s^2,\mu_\star\},
% \qquad
% \epsilon_s=\frac{\gamma}{\mu_s}.
% \end{equation}
% Let
% \[
% \gamma=\frac{\sqrt{40(M_f^2+M_h^2+m^2\rho^2\bar{M}_c^2M_c^2)}}{\beta},
% \qquad
% \mu_{\rm base}=\max\{1,\gamma,
% \sqrt{\frac{10(8L_v^2+3)}{27\beta^2}},
% \sqrt{\frac{8(2L_v^2+1)}{3\beta^2}}\}.
% \] 
%At stage \(s\), the stepsizes for dual update are \(\rho_{s,k}=\rho/K_s\), \(k=0,\ldots,K_s-1\). 
% The parameters \(\nu_s\), \(\eta_s\), \(\alpha_s\), \(n_{0,s}\), and \(n_s\) are chosen according to the corresponding parameter setting in Theorem~\ref{thm:p} (or Theorem~\ref{thm:pinf}), with \(\epsilon\) replaced by \(\epsilon_s\) and \(\mu\) replaced by \(\mu_s\). 
By the algorithm framework, it is easy to obtain that, if \(\mu_\star=\mu_1=\mu_{\rm base}\), then \(T=1\). Otherwise, before the algorithm terminates, it always happens that \(\mu_{s+1}=2\mu_s^2\). It derives \(\mu_s=\frac12(2\mu_1)^{2^{s-1}}\), thus
\be\label{O_T}
T\le1+\lceil\log_2(\frac{\log(2\mu_\star)}{\log(2\mu_1)})\rceil
=\mathcal O(1+\log\log\frac{e}{\epsilon}).\ee
Moreover, at each  stage all samples are freshly generated independent and identically distributed  samples according to Algorithm~\ref{algo}; conditioned on the history before the current stage, they follow the same sampling rule as in the single-stage analysis.

Let the iteration number $K_s$ in $s$-th stage and the according dual update stepsizes $\rho_s$ satisfy
\be\label{K-s}
K_s=\lceil\max\{\frac{108\bar\Delta_s}{\eta_s\epsilon_s^2},
\frac{5\bar C_s}{\epsilon_s^2},
\frac{1}{L_f^2\eta_s^2}\}\rceil,
\quad \rho_{s,k}=\frac{\rho}{K_s},\quad \rho>0,\ k=0,\ldots,K_s-1.
\ee
with
\be\label{par-set-s}
\begin{gathered}
\gamma=\frac{\sqrt{40(M_f^2+M_h^2+m^2\rho^2\bar{M}_c^2M_c^2)}}{\beta},\quad \mu_{\rm base}=\max\{1,\gamma,
\sqrt{\frac{10(8L_v^2+3)}{27\beta^2}},
\sqrt{\frac{8(2L_v^2+1)}{3\beta^2}}\},\\
B_X=\sup_{x\in X}\{f(x)+h(x)-\Phi_{\rm inf}\}+2m\rho \bar{M}_c^2,\ 
 \bar\Delta_1=B_X+\mu_1m\bar{M}_c^2/2, \ \bar\Delta_+=B_X+m^{1-2/p}\gamma^2,\\
\bar\Delta_s=
\begin{cases}
\bar\Delta_1, & s=1,\\[1mm]
\bar\Delta_+, & s\ge2,
\end{cases}
\qquad
\bar C_1=m\bar{M}_c^2,
\qquad
\bar C_s=\frac{2m^{1-2/p}\gamma^2}{\mu_s},\quad s\ge2.
\end{gathered} 
\ee
By the algorithm framework and  above parameter settings, 
\(\mu_1=\mu_{\rm base}=\mathcal O(1+L_v)\) and \(\bar\Delta_1=\mathcal O(1+L_v)\), whereas for \(s\ge2\), \(\bar\Delta_s=\bar\Delta_+=B_X+m^{1-2/p}\gamma^2=\mathcal O(1)\). Thus, only the first stage carries a possible \(L_v\)-dependent initial augmented-Lagrangian gap.

As will be shown in Theorem \ref{thm:restart}, without the near initial feasibility the leading {high-accuracy} oracle complexity same as the single-stage algorithm is still  preserved. The additional cost caused by an arbitrary initial point appears only through an \(\epsilon\)-independent first-stage term. The proof  is deferred to Appendix~\ref{app:restart}. 

\begin{theorem}
\label{thm:restart}
 Let every stage of Algorithm~\ref{algo:restart} use the same mirror descent metric generated by \(v_\delta\) in \eqref{eq:smoothed_v}.  Suppose that \(u\) follows the Rademacher distribution,  Assumptions \ref{ass:lip}--\ref{ass:ms-p} hold, and  Assumption~\ref{ass:ensc} holds at every iterate generated in every stage, with a common constant \(\beta>0\) across all stages. Further assume that the bound $D_X$ in Assumption~\ref{ass:lip} and the regularity constant $\beta$ in Assumption~\ref{ass:ensc} are  independent of \(d\) and \(p\).  Under parameter settings \eqref{K-s} and \eqref{par-set-s},  at $s$-th stage the parameters \(\nu_s\), \(\eta_s\), \(\alpha_s\), \(n_{0,s}\), \(n_s\) are chosen following the corresponding parameter settings in Theorem~\ref{thm:p}  when \(p\in[2,2\ln d]\) and in Theorem~\ref{thm:pinf} when \(p>2\ln d\), with \(\epsilon\) replaced by \(\epsilon_s\) and \(\mu\) replaced by \(\mu_s\).   Then, Algorithm \ref{algo:restart} returns an \(\epsilon\)-KKT point %after \(T=\mathcal O(1+\log\log(e/\epsilon))\) stages
with oracle complexity in order
\be\label{mul-cp}
\begin{cases}
\mathcal O(p d^{2/p}\epsilon^{-3}+p^5d^{4-6/p}), & p\in[2,2\ln d],\\[2mm]
\mathcal O(\ln d\cdot\epsilon^{-3}+p^4d^{4-8/p}\ln d), & p>2\ln d.
\end{cases}
\ee
\end{theorem}

In the Euclidean case, the additional first-stage initialization cost is absorbed into the leading {high-accuracy} oracle complexity term, yielding the same oracle complexity as in Corollary~\ref{thm:2-norm}. We also present the proof of Corollary \ref{cor:2-restart}  in Appendix \ref{app:restart}.

\begin{corollary}\label{cor:2-restart}
Let every stage of Algorithm~\ref{algo:restart} use the same mirror descent metric generated by \(v(x) = \frac{1}{2}\|x\|_2^2\). Suppose that \(u\) follows the Rademacher distribution,  Assumptions \ref{ass:lip}--\ref{ass:ms-p} hold, and  \eqref{CQ_2} holds at every iterate generated in every stage, with a common constant \(\beta>0\) across all stages. Under parameter settings \eqref{K-s} and \eqref{par-set-s} with \(p=2, L_v=1\), at $s$-th stage the parameters \(\nu_s\), \(\eta_s\), \(\alpha_s\), \(n_{0,s}\), \(n_s\) are chosen following the corresponding parameter settings in Corollary \ref{thm:2-norm} with \(\epsilon\) replaced by \(\epsilon_s\) and \(\mu\) replaced by \(\mu_s\). Algorithm \ref{algo:restart} returns an \(\epsilon\)-KKT point 
%after \(T=\mathcal O(1+\log\log(e/\epsilon))\) stages, 
with oracle complexity in order \(\mathcal O(d\epsilon^{-3})\). 
\end{corollary}

\section{Numerical experiments}\label{sec:num}
In this section, we report three sets of numerical experiments. The first experiment uses a controlled QCQP dimension-dependence ablation to verify that the non-Euclidean geometry reduces the high-dimensional query dependence predicted by the theory. The second experiment tests whether this query advantage persists in a more realistic MNIST black-box adversarial attack. The third experiment evaluates the Euclidean version of the proposed method on a nonconvex fairness-constrained classification problem and compares it with existing stochastic zeroth-order baselines. In the first two experiments, all variants, including the non-Euclidean one, are evaluated using the Euclidean 
$\ell_2$ KKT metric.

\subsection{Dimension-dependence ablation}
\label{subsec:dimension-ablation}

We first isolate the dimension dependence predicted by {Theorem~\ref{thm:rad_qnorm}} on a controlled full-dimensional
synthetic problem. Unlike constructions that artificially increase the dimension by adding inactive coordinates, we use a dense indefinite quadratic objective in which all coordinates are actively involved.

For each
\(d\in\{2^6,2^8,2^{10},2^{12},2^{14}\}\), let \(H_d\in\mathbb R^{d\times d}\) be the normalized
Walsh--Hadamard matrix, so that \(H_d^\top H_d=I_d\) and every entry of \(H_d\) is \(\pm1/\sqrt d\).
Let \(\Lambda_d\) be a diagonal matrix whose first \(d/2\) prescribed eigenvalues are evenly spaced in
\([-0.05,-0.01]\) and whose remaining \(d/2\) eigenvalues are evenly spaced in \([0.01,0.05]\), followed by a
fixed random permutation. For each random seed, we   
define
\(
%\beta_d=\frac{\tilde{\beta}_d}{\|\tilde{\beta}_d\|_2}, 
A_d=H_d\Lambda_dH_d^\top, 
X_d=\{x\in\mathbb R^d:\|x\|_2\le2.5\}, 
b_d=H_d \tilde{\beta}_d/\|\tilde{\beta}_d\|_2
\) with $\tilde{\beta}_d$ following standard normal distribution $\mathcal{N}(0,I_d)$. 
Generate the test problem 
\[
\min_{x\in X_d}\quad f_d(x):=\frac12x^\top A_dx+b_d^\top x
\qquad\text{s.t.}\qquad
c_d(x):=\frac12\bigl(\|x\|_2^2-1\bigr)=0.
\]
The compact convex set \(X_d\) matches the boundedness assumption in the theory. Trial iterates with norm
larger than \(2.5\) are projected onto its boundary. The deterministic oracle is the special case
\(F_d(x;\xi)\equiv f_d(x)\). 
%On the feasible manifold, \(\|x\|_2=1\) and \(\nabla c_d(x)=x\ne0\), so linear independence constraint qualification (LICQ) holds; the corresponding non-singularity condition holds locally in neighborhoods of the sphere that stay away from the origin. 
All methods start from the same feasible point
$x^0=\mathbf{1}/\sqrt{d}$ and use $\lambda^0=0$.

We compare the following three variants of Algorithm~\ref{algo}.
\begin{itemize}
    \item \textbf{Euc-theory:} the Euclidean geometry with $p=q=2$,
    $v(x)=\frac{1}{2}\|x\|_2^2$, and $n_0=n=d$ random directions per iteration;

    \item \textbf{Euc-matched:} the same Euclidean geometry, but with the number of directions matched
    to that of the non-Euclidean method;

    \item \textbf{Non-Euc:} the non-Euclidean geometry with \(p_d=2\ln d,\ q_d=\frac{p_d}{p_d-1}\),
    and the smoothed distance-generating function \eqref{eq:smoothed_v} by \(q_d,\ \delta=10^{-2}\).
\end{itemize}

At every iteration, all methods independently sample \( u_j^k\in\{-1,1\}^d,\ j=1,\ldots,n,\)
from the Rademacher distribution. The directions are mutually independent across both $j$ and $k$.
Consequently, \(\mathbb{E}[u_j^k(u_j^k)^\top]=I_d,\ \mathbb{E}
[\|\langle g,u_j^k\rangle u_j^k\|_p^2]=d^{2/p}\|g\|_2^2.\)
Thus $S_2=d$ in the Euclidean geometry, whereas
$S_{p_d}=d^{2/p_d}=e$. In the non-Euclidean geometry, %following the batch-size in \eqref{settings},
we set \(n_0=n
=\lceil(p_d-1)d^{2/p_d}\rceil=\lceil e(p_d-1)\rceil\). Both Euc-matched and Non-Euc  use exactly
the same number of random directions and function evaluations per iteration.
The same stepsize schedule \(\eta_k=\frac{\eta_0}{\sqrt{1+k/25}}\) is used for all three methods. The base stepsize $\eta_0$ is selected from the same prescribed grid \(\{0.005,0.01,0.02,0.04,0.08\}\)
using independent pilot runs, and the selected values are then fixed for all reported random seeds.
The remaining parameters are $\mu=20$, $\nu=10^{-6}$, and $\alpha=0.5$.

Although the non-Euclidean variant is motivated by an $\ell_p$ KKT metric, to compare the methods under a common criterion, we use the same Euclidean KKT diagnostic for all
three variants. Define \(g_{\rm AL}^k:=\nabla f_d(x^k)+(\lambda^k+\mu c_d(x^k))\nabla c_d(x^k),\) and the Euclidean projected-gradient mapping \(\mathcal{G}_{X_d}^{(2)}(x,g):=x-\Pi_{X_d}(x-g).\)
The reported residual is \(r_2^k:=\max\{\|\mathcal{G}_{X_d}^{(2)}(x^k,g_{\rm AL}^k)\|_2,|c_d(x^k)|\}.\)
For every run, we record the first function-query count at which $r_2^k\le 10^{-2}$. We use ten paired random seeds, with the same problem instance and initial
point supplied to all methods under each seed.

\begin{figure}[ht]
\centering
\includegraphics[width=.72\linewidth]{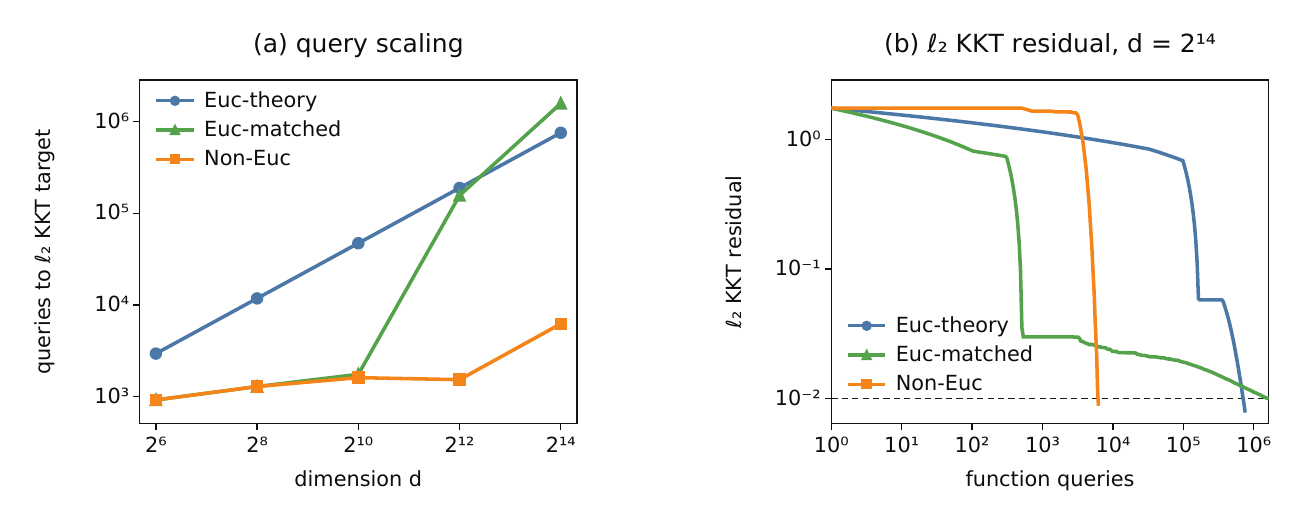}
\caption{QCQP dimension-dependence ablation. (a) Queries required to reach the Euclidean
\(\ell_2\) KKT target as the dimension increases. (b) \(\ell_2\) KKT residual versus function
queries at \(d=2^{14}\). }
\label{fig:qcqp-kkt-scaling}
\end{figure}

Figure~\ref{fig:qcqp-kkt-scaling}(a) reports the mean number of function queries required to reach Euclidean $\ell_2$ KKT target. Both Euc-theory and Euc-matched exhibit pronounced dimension dependence, whereas the Non-Euc  remains substantially more query-efficient in high dimensions. At $d=2^{14}$, Euc-matched  requires $1.58\times10^6$ queries, whereas Non-Euc requires $6.22\times10^3$ queries,  approximately $253$ times fewer queries.
Figure~\ref{fig:qcqp-kkt-scaling}(b) shows three methods corresponding \(\ell_2\) KKT residual decrease at
\(d=2^{14}\).  It plots the best residual attained at each query budget. And the curves show medians over ten paired seeds. For  Euc-theory, its use of $d$ directions per iteration leads to a high query cost, and it reaches the target after approximately $7.54\times10^5$ queries. Euc-matched achieves a relatively small residual early but then remains above the target for a long query interval, requiring approximately $1.58\times10^6$ queries to terminate. In contrast, after a short initial phase, Non-Euc exhibits a rapid residual decrease and reaches the target with only $6.22\times10^3$ queries. 

\subsection{Black-box adversarial attack}
\label{subsec:black-box-attack}

We next test whether the query advantage of the non-Euclidean geometry appears in a more realistic
black-box task. Since the MNIST dimension is fixed at \(d=784\), this experiment serves as a
complementary high-dimensional validation. The problem
corresponds to the special case of Algorithm~\ref{algo} with no equality constraint, \(c\equiv 0\), and with the
box constraint handled by the mirror subproblem.

We train a softmax classifier on 20,000 MNIST training images, and let \(\mathcal{M}_j(x)\) denote
the logit of class \(j\) at input \(x\). The resulting clean test accuracy is \(90.76\%\). We attack 40
correctly classified test images with clean margin
\(\mathcal{M}_y(x_0)-\max_{j\ne y}\mathcal{M}_j(x_0)>0.2\). Given an image \(x_0\in[0,1]^d\)
with label \(y\), the feasible set is
\(
    X(x_0)=\{\delta:\|\delta\|_\infty\le 0.3,\;
    0\le x_0+\delta\le 1\}.
\)
The optimizer minimizes the smoothed margin
\[
f_\tau(\delta)
=
\mathcal{M}_y(x_0+\delta)
-
\tau\log\sum_{j\ne y}
\exp(\frac{\mathcal{M}_j(x_0+\delta)}{\tau}),
\]
with \(\tau=0.5\), while attack success is evaluated using the nonsmoothed classification margin
\[
    \mathcal{M}_y(x_0+\delta)-\max_{j\ne y}\mathcal{M}_j(x_0+\delta)<0.
\]
The optimizer has access only to function values of \(f_\tau\). For consistency with the ablation,
the non-Euclidean variant is evaluated under the same Euclidean diagnostic. And the reported attack metrics are the common model-query budget,
success indicator, and nonsmoothed classification margin. All variants use the same query budget of
12,000 model queries and momentum parameter \(0.5\).
We report three methods: the Euc-theory variant with \(p=2\) and \(n_0=n=d=784\) directions,
the Euc-matched  variant with \(p=2\) and \(n_0=n=69\) directions, and the Non-Euc variant
with \(p=\lceil 2\ln d\rceil+1=15\) and
\(n_0=n=\lceil 2(p_d-1)d^{2/p_d}\rceil=69\) directions.

\begin{figure}[t]
    \centering
    \includegraphics[width=.72\linewidth]{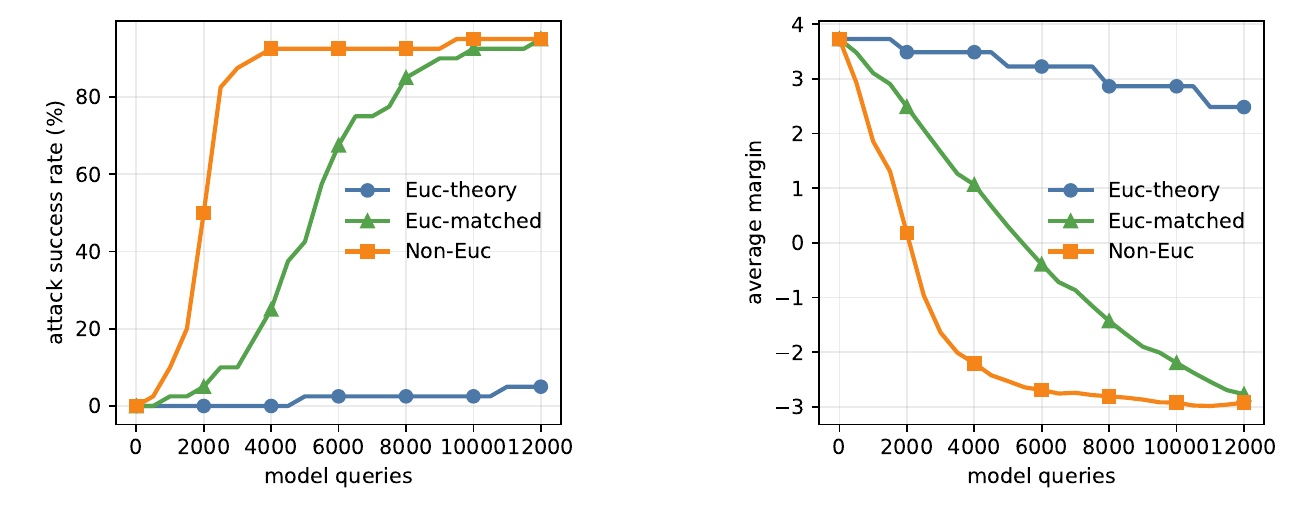}
    \caption{Black-box untargeted attack on MNIST. Left: attack success rate as a function of model
    queries. Right: average nonsmoothed classification margin as a function of model queries.}
    \label{fig:mnist-black-box-attack}
\end{figure}

\begin{table}[ht]
\centering
\caption{Black-box untargeted attack on MNIST. ``Dirs./it." denotes random directions per iteration. ``Succ." reports the number of successful attacks within
the query budget.  ``Median queries (Succ.)" is the median number of model queries among successful attacks only.  ``Median queries (all
images)" is the censored median query count over all 40 images, where failed attacks are assigned the
full 12,000-query budget. ``Final margin" is the final true-class logit minus the largest competing
logit, reported as mean \(\pm\) standard deviation over all images. }
\label{tab:mnist-black-box-attack}
\begingroup
\small
\setlength{\tabcolsep}{2.5pt}
\newcommand{\tableheadcell}[1]{%
    \begin{tabular}[c]{@{}c@{}}#1\end{tabular}}
\begin{tabular}{@{}l
    S[table-format=2.0]
    S[table-format=3.0]
    S[table-format=2.0]@{/}S[table-format=2.0]
    S[table-format=5.0]
    S[table-format=5.0]
    S[table-format=-1.3]@{\,\(\pm\)\,}S[table-format=1.3]@{}}
\toprule
\multicolumn{1}{c}{\tableheadcell{Method}} &
\multicolumn{1}{c}{\tableheadcell{$p$}} &
\multicolumn{1}{c}{\tableheadcell{Dirs./it.}} &
\multicolumn{2}{c}{\tableheadcell{Succ.}} &
\multicolumn{1}{c}{\tableheadcell{Median queries\\(Succ.)}} &
\multicolumn{1}{c}{\tableheadcell{Median queries\\(all images)}} &
\multicolumn{2}{c}{\tableheadcell{Final margin}}\\
\midrule
Euc-theory  &  2 & 784 &  2 & 40 & 7840 & 12000 &  2.484 & 1.736 \\
Euc-matched &  2 &  69 & 38 & 40 & 5382 &  5382 & -2.770 & 1.520 \\
Non-Euc    & 15 &  69 & 38 & 40 & 1794 &  1932 & -2.924 & 1.667 \\
\bottomrule
\end{tabular}
\endgroup
\end{table}

Figure~\ref{fig:mnist-black-box-attack} shows how attack performance evolves with  the query budget. In the left penal, Non-Euc reaches  \(50\%\) attack success rate after \(2000\) quaries and
\(92.5\%\) after \(4000\) queries. By comparison, Euc-matched attains only \(5\%\) and \(25\%\) attack success at the
same query budgets and increases more gradually to \(95\%\) at \(12{,}000\) queries, while Euc-theory 
remains at or below \(5\%\) throughout the experiment. In the right panel, the average margin of Non-Euc becomes negative shortly after \(2000\) queries and is already close to \(-2.2\) at \(4000\) queries. The average margin of Euc-matched becomes negative near \(6000\) queries, whereas that of Euc-theory remains positive over the entire query budget. Table \ref{tab:mnist-black-box-attack} quantifies these observations at termination. Euc-theory succeeds on only \(2\) of the \(40\) images. Euc-matched and Non-Euc both succeed on \(38\) images and
use the same \(69\) random directions per iteration, but their median query counts among successful
attacks are \(5382\) and \(1794\), respectively. Thus, Non-Euc achieves the same success count with
approximately three times fewer median queries.  Its median query count over all images is also lower, and its final average margin is slightly more negative. Hence, the query reduction is not obtained by sacrificing attack success
or final attack strength, and it cannot be explained solely by using fewer directions per iteration, but rather by the non‑Euclidean mirror descent geometry used in the Non‑Euc, which better captures the geometric structure of the test  problem. 

\subsection{Nonconvex fairness constrained problem}
Finally, we evaluate the Euclidean version of Algorithm \ref{algo} on the nonconvex fairness-constrained problem used by Li et al.~\cite{li2024stochastic}. Write the labeled training set as
\(\mathcal D=\{(a_i,b_i)\}_{i=1}^{N}\), where \(a_i\in\mathbb R^d\) and \(b_i\in\{-1,1\}\), and define
\[
   f_0(x)=\frac{1}{N}\sum_{i=1}^{N}
2\log(1+\frac{1}{2}\log(1+\exp(-b_i a_i^\top x))).
\]
Let \(S\) be the samples on which fairness is evaluated and let
\(S_{\min}\subseteq S\) denote the minority-group subset. With a prescribed fairness weight
\(\omega>0\), set
\[
    g(x)=
\frac{\omega}{|S|}\sum_{a\in S}\sigma(a^\top x)
-\frac{1}{|S|}\sum_{a\in S_{\min}}\sigma(a^\top x),
\qquad
\sigma(t)=\frac{\exp(t)}{1+\exp(t)}.
\]
The original fairness requirement is \(g(x)\le0\). To express it exactly in the equality-constrained form \eqref{p}, introduce a slack variable \(s\ge0\), let \(z=(x,s)\) and \(X_f=\mathbb R^d\times\mathbb R_+\), and solve
\[
\min_{z\in X_f}\ f(z):=f_0(x)
\qquad\text{s.t.}\qquad
 c(z):=g(x)+s=0.
\]
The dimensions and sample counts in the
three released data sets are
\[
\begin{array}{c|rrrr}
\text{data set} & d & N & |S| & | S_{\min}|\\ \hline
\textit{a9a}  & 123 & 32561 & 16281 & 1561\\
\textit{bank} & 81 & 22605 & 22605 & 6320\\
\textit{loan} & 250 & 64485 & 63890 & 31966
\end{array}
\]

We apply Algorithm~\ref{algo} to the equality-constrained variable
\(z=(x,s)\in X_f\), retaining and updating the slack variable throughout the run.
To avoid conflict between the slack variable \(s\) and the stochastic gradient estimator
denoted by \(s^k\) in Algorithm~\ref{algo}, we write the estimator in this paragraph as
\(v^k=(v_x^k,0)\). Since \(f(z)=f_0(x)\) is independent of the slack variable, the slack
component of the stochastic gradient estimator is zero; the \(x\)-component \(v_x^k\)
is computed by the recursion \eqref{hsgd}, using the two-function-value Rademacher
estimator \eqref{sgrad}, \(8\)  independent random directions, and a mini-batch of
\(256\) data samples. 
We set \(p=2\) and use the Euclidean Bregman divergence on the full variable,
\[
V(z,w)=\frac12\|z-w\|_2^2
=\frac12\|x-y\|_2^2+\frac12(s-t)^2,
\qquad z=(x,s),\quad w=(y,t).
\]
Due to 
\(
\nabla c(z^k)=(\nabla g(x^k),1),
\)
the primal and dual update have the explicit expression:
\[
\begin{gathered}
\theta^k=\lambda^k+\mu\bigl(g(x^k)+s^k\bigr),\quad 
x^{k+1}
=x^k-\eta(v_x^k+\theta^k\nabla g(x^k)),\quad 
s^{k+1}
=[s^k-\eta\theta^k]_+,\\
\lambda^{k+1}
=\lambda^k+\rho_k c(z^k)
=\lambda^k+\rho_k\bigl(g(x^k)+s^k\bigr).
\end{gathered}
\]
We set \(\nu=10^{-3}\), parameter \(\alpha=1\), penalty parameter
\(\mu=1\), and \(\rho_k=0.5\). Every run is initialized with
\(
x^0=0, s^0=[-g(x^0)]_+, \lambda^0=0.
\)
The primal stepsize \(\eta\) is held fixed within each run and selected by the
same validation protocol; the values are \(0.60\), \(2.00\), and \(0.40\) for
\textit{a9a}, \textit{bank}, and \textit{loan}, respectively.
The proposed method is terminated and
diagnosed by  residuals
\(
\mathrm{pres}(z)=|c(z)|=|g(x)+s|
\)
and the generalized gradient residual
\(
\mathrm{dres}(z)
=\|
\mathcal G_{X_f}(
z,\,
(\nabla f_0(x),0)
+\tilde\lambda\,(\nabla g(x),1),\,
\eta)\|_2\) with \(
\tilde\lambda=\lambda+\mu c(z).
\)

For comparison, we implement the Euclidean specializations of SZO-ConEx~\cite{nguyen2022stochastic}
and the stochastic zeroth-order penalty method ZO-Penalty~\cite{wang2017penalty} with the same zeroth-order oracle,
smoothing parameter, and mini-batch size. All methods start from $x_0=0$, and each result is averaged over
ten paired random seeds. Following the stopping criteria in \cite{li2024stochastic}, we report the feasibility, stationarity and complementary slackness: 
\[
\mathrm{Feas.}=[g(x)]_+,\quad
\mathrm{Stat.}=\|\nabla f_0(x)+\widehat\lambda(x)\nabla g(x)\|_2,\quad
\mathrm{Cs.}=|\widehat\lambda(x)g(x)|,
\]
where
\(
\widehat\lambda(x)=
 [-\frac{\langle\nabla f_0(x),\nabla g(x)\rangle}
{g(x)^2+\|\nabla g(x)\|_2^2} ]_+ .
\)

\begin{figure}[ht]
    \centering
    \includegraphics[width=.95\linewidth]{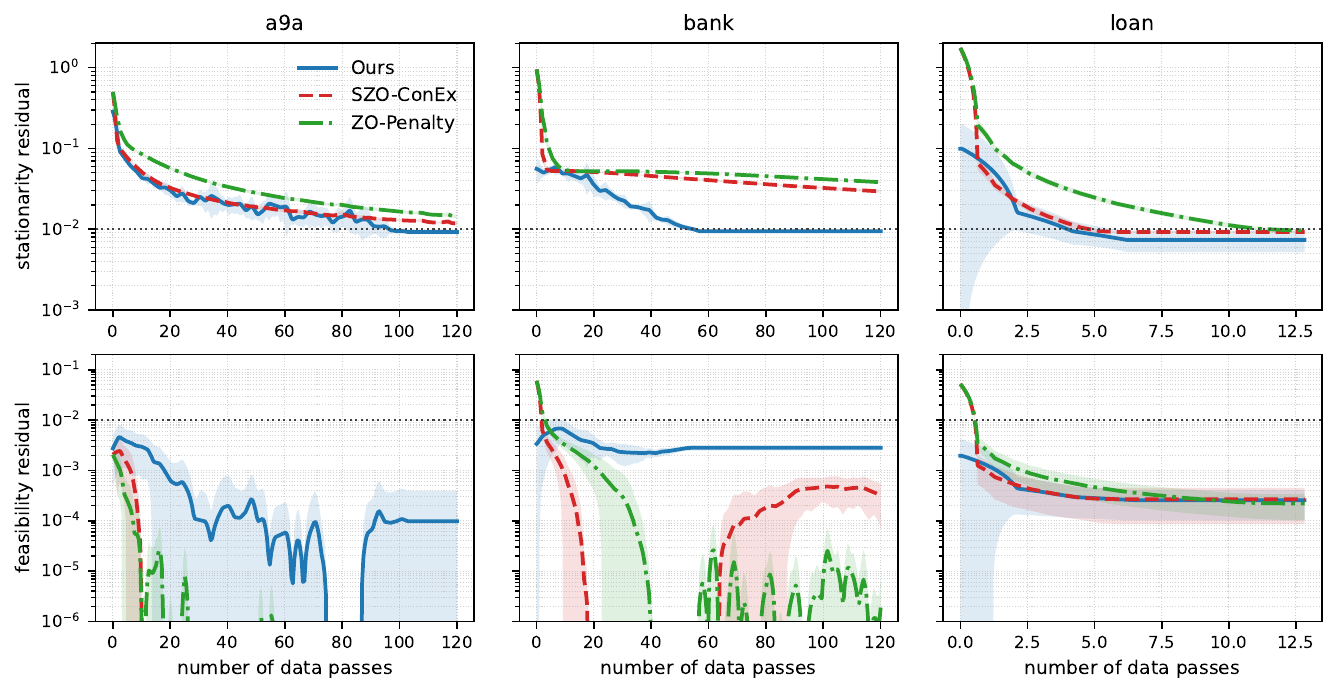}
    \caption{Fairness constrained comparison. Top: common stationarity residuals. Bottom:
    common feasibility residuals. Columns correspond to \textit{a9a}, \textit{bank}, and
    \textit{loan}.}
    \label{fig:fairness-comparison}
\end{figure}

\begin{table}[ht]
\centering
\caption{Final common KKT diagnostics and data passes on the fairness constrained problem.
All results are reported as mean \(\pm\) standard deviation over ten paired random seeds.}
\label{tab:fairness-zo-results}
\begingroup
\small
\setlength{\tabcolsep}{2.5pt}
\begin{tabular}{@{}ll
    r@{\,\(\pm\)\,}l
    r@{\,\(\pm\)\,}l
    r@{\,\(\pm\)\,}l
    r@{\,\(\pm\)\,}l@{}}
\toprule
Data & Method & \multicolumn{2}{c}{Feas.} & \multicolumn{2}{c}{Stat.} &
\multicolumn{2}{c}{Cs.} & \multicolumn{2}{c}{Data passes}\\
\midrule
\textit{a9a}  & Ours       & 9.76e-05 & 3.09e-04 & 9.26e-03 & 6.44e-04 & 3.36e-05 & 5.46e-05 &  85.46 & 13.82\\
\textit{a9a}  & SZO-ConEx  & 0.00e+00 & 0.00e+00 & 1.19e-02 & 6.86e-04 & 2.47e-05 & 3.35e-05 & 120.01 &  0.00\\
\textit{a9a}  & ZO-Penalty & 0.00e+00 & 0.00e+00 & 1.45e-02 & 7.43e-04 & 2.17e-06 & 6.51e-06 & 120.01 &  0.00\\
\textit{bank} & Ours       & 2.82e-03 & 3.55e-04 & 9.46e-03 & 3.05e-04 & 1.03e-03 & 2.31e-04 &  54.54 &  2.51\\
\textit{bank} & SZO-ConEx  & 3.17e-04 & 2.46e-04 & 2.94e-02 & 9.12e-04 & 1.04e-04 & 7.22e-05 & 120.13 &  0.00\\
\textit{bank} & ZO-Penalty & 1.98e-06 & 5.95e-06 & 3.83e-02 & 7.24e-04 & 1.73e-05 & 9.73e-06 & 120.13 &  0.00\\
\textit{loan} & Ours       & 2.56e-04 & 1.64e-04 & 7.36e-03 & 2.27e-03 & 9.06e-06 & 2.86e-05 &   4.78 &  1.95\\
\textit{loan} & SZO-ConEx  & 2.68e-04 & 1.82e-04 & 9.25e-03 & 1.82e-04 & 5.66e-06 & 1.70e-05 &   5.06 &  0.45\\
\textit{loan} & ZO-Penalty & 2.20e-04 & 1.17e-04 & 9.69e-03 & 2.10e-04 & 5.23e-08 & 1.57e-07 &  11.41 &  0.75\\
\bottomrule
\end{tabular}
\endgroup
\end{table}
Figure~\ref{fig:fairness-comparison} compares the stationarity and feasibility residual curves of the proposed method with those of SZO-ConEx~\cite{nguyen2022stochastic} and ZO-Penalty~\cite{wang2017penalty}. In the top row, the proposed method reaches the $10^{-2}$ stationarity threshold earlier than the two baselines on \textit{a9a} and \textit{bank}; on \textit{bank}, the baseline residuals remain above this threshold even after approximately $120$ data passes. On \textit{loan}, all three methods reach residuals near $10^{-2}$, but the proposed method achieves this with the fewest data passes.  The bottom row shows that all methods maintain feasibility residuals below the $10^{-2}$ threshold. Although the baselines attain smaller feasibility residuals on \textit{a9a} and \textit{bank}, this comes with slower stationarity reduction. Table~\ref{tab:fairness-zo-results} quantifies these observations. The proposed method achieves the smallest stationarity residual on all three data sets and uses fewer data passes than the two baselines. All three methods obtain small feasibility and complementary-slackness residuals, although the baselines are smaller in some cases. %Overall, the proposed method reaches better stationarity with less computation.

\section{Conclusions}
\label{sec:conclusions}
In this paper, we study a class of nonconvex constrained stochastic optimization problems for which exact constraint information is available but the objective can be accessed only through noisy function values. We propose a single-loop mirror descent linearized augmented Lagrangian method. At each iteration, a momentum-based variance-reduction recursion is incorporated into  a two-point stochastic zeroth-order estimator to approximate the objective gradient. The smooth part of the augmented Lagrangian is then linearized, and the primal variable is updated by solving a Bregman proximal subproblem. We establish oracle complexity guarantees for finding an $\epsilon$-KKT point under general mirror descent geometries. In particular, under Rademacher smoothing, the leading high-accuracy oracle complexity order is 
\(\mathcal{O}(p d^{2/p}\epsilon^{-3})\) for \(p\in[2,2\ln d]\) and
\(\mathcal{O}(\ln d\cdot\epsilon^{-3})\) for \(p>2\ln d\).
In the standard Euclidean setting with \(p=2\), the proposed method achieves an oracle complexity of
\(\mathcal{O}(d\epsilon^{-3})\). We also present a multi-stage scheme to remove the near‑feasibility of the initial point while maintaining the leading high-accuracy oracle complexity orders. 
These results show that an appropriate non-Euclidean mirror geometry can potentially reduce the dimension dependence of zeroth-order gradient estimation. Finally, three sets of numerical experiments demonstrate the effectiveness and query efficiency of the proposed method.

\bibliographystyle{siamplain}
\bibliography{references}
\appendix
\section{Properties of the smoothed mirror map}
\label{app:smooth-mirror}
\begin{lemma}\label{lem:smooth-mirror-map}
Let $v_\delta$ be defined by \eqref{eq:smoothed_v}. %, where $\frac1p+\frac1q=1$.
Then $v_\delta$ satisfies Assumption \ref{ass:v}
%is 1-strongly convex w.r.t. the $\ell_q$-norm, and $\nabla v_\delta$ is $L_{v}$-Lipschitz continuous on $X$, where
with $L_v$ defined in \eqref{eq:Lv}.
\end{lemma}
\begin{proof}
By Assumption \ref{ass:lip}, $X$ is bounded in the $\ell_q$-norm with $D_X$.
Let
\( \Psi_\delta(x)=\sum_{i=1}^d (x_i^2+\delta^2)^{q/2},\) and \(a_i=x_i^2+\delta^2 .\)
Then, we obtain
\(
    v_\delta(x)=\frac{1}{2(q-1)}\Psi_\delta(x)^{2/q}.
\) 
A direct calculation gives the gradient entries
\[
    \nabla_i v_\delta(x)=\frac{1}{q-1}\Psi_\delta(x)^{2/q-1}a_i^{q/2-1}x_i ,\quad i=1,\ldots,d.
\]
Moreover, the Hessian entries are
\[
    [\nabla^2 v_\delta(x)]_{ij}
    =\frac{2-q}{q-1}\Psi_\delta(x)^{2/q-2}a_i^{q/2-1}x_ia_j^{q/2-1}x_j                +\frac{1}{q-1}\Psi_\delta(x)^{2/q-1}a_i^{q/2-2}\big((q-1)x_i^2+\delta^2\big)\mathbf 1_{{i=j}} 
\]
for $i,j=1,\ldots,d.$ 
Thus, we obtain\(\nabla^2 v_\delta(x)=M(x)+\Lambda(x),\) where \(M(x)=\frac{2-q}{q-1}\Psi_\delta(x)^{2/q-2}ww^\top, w_i=a_i^{q/2-1}x_i,\) and \(\Lambda_{ii}(x)=\frac{1}{q-1}\Psi_\delta(x)^{2/q-1}a_i^{q/2-2}((q-1)x_i^2+\delta^2), \Lambda_{ij}(x)=0, \) for \(i,j=1,\ldots,d;\ i\ne j .\)

We first prove 1-strong convexity. Since $M(x)$ is positive semidefinite, for any $y\in\mathbb R^d$,
\[
    y^\top \nabla^2 v_\delta(x)y
    \ge
    y^\top \Lambda(x)y     =
    \frac{1}{q-1}
\Psi_\delta(x)^{(2-q)/q}
    \sum_{i=1}^d
    a_i^{q/2-2}
    \big((q-1)x_i^2+\delta^2\big)y_i^2    
\ge\Psi_\delta(x)^{(2-q)/q}
    \sum_{i=1}^d
    a_i^{(q-2)/2}y_i^2 .
\]
The last inequality follows from\(\frac{(q-1)x_i^2+\delta^2}{q-1}=x_i^2+\frac{\delta^2}{q-1}\ge x_i^2+\delta^2=a_i .\) On the other hand, by H\"older's inequality, it holds
\[
    \sum_{i=1}^d |y_i|^q
    =
    \sum_{i=1}^d
    (a_i^{(q-2)/2}y_i^2)^{q/2}
    (a_i^{q/2})^{(2-q)/2}       
    \le
    (\sum_{i=1}^d a_i^{(q-2)/2}y_i^2)^{q/2}
    (\sum_{i=1}^d a_i^{q/2})^{(2-q)/2}.
\]
Increase both sides to the power $2/q$ yields \(\|y\|_q^2\le\Psi_\delta(x)^{(2-q)/q}\sum_{i=1}^d a_i^{(q-2)/2}y_i^2 .\) Consequently, we have \(y^\top \nabla^2 v_\delta(x)y\ge\|y\|_q^2 ,\ \forall x,y\in\mathbb R^d\), applying the fundamental theorem of calculus gives
\(\langle x-y,\nabla v_\delta(x)-\nabla v_\delta(y)\rangle=\int_0^1(x-y)^\top\nabla^2v_\delta(y+t(x-y))(x-y)\,dt\geq
\|x-y\|_q^2\), which proves that $v_\delta$ is 1-strongly convex w.r.t. the $\ell_q$-norm.

We next proof $v_\delta$ is \(L_v\)-smooth on $X$.
 H\"older's inequality gives
\(
    \langle w,h\rangle^2
    \le
    \|w\|_p^2\|h\|_q^2.
\) 
Since \(|w_i|=a_i^{q/2-1}|x_i|\le a_i^{(q-1)/2}\), we have \(\|w\|_p^p\le\sum_{i=1}^d a_i^{p(q-1)/2}=\sum_{i=1}^d a_i^{q/2}=\Psi_\delta(x),\) where $p=q/(q-1)$. Hence \(\|w\|_p^2\le\Psi_\delta(x)^{2/p}=\Psi_\delta(x)^{2(q-1)/q}.\) Therefore,
\[
    y^\top M(x)y
    =
    \frac{2-q}{q-1}
 \Psi_\delta(x)^{2/q-2}
    \langle w,y\rangle^2  
    \le
    \frac{2-q}{q-1}   \Psi_\delta(x)^{2/q-2}    \Psi_\delta(x)^{2(q-1)/q}
    \|y\|_q^2        =
    \frac{2-q}{q-1}\|y\|_q^2 .
\]
On the other hand, since $q\in(1,2]$, we have \(\frac{(q-1)x_i^2+\delta^2}{x_i^2+\delta^2}\le 1,\) and \(a_i^{(q-2)/2}\le\delta^{q-2},\) thus
\[
    \Lambda_{ii}(x)
    =
    \frac{1}{q-1}
    \Psi_\delta(x)^{(2-q)/q}
    a_i^{(q-2)/2}
    \frac{(q-1)x_i^2+\delta^2}{x_i^2+\delta^2}              \le
    \frac{1}{q-1}
    \Psi_\delta(x)^{(2-q)/q}
    \delta^{q-2}.
\]
Moreover, for $x\in X$, we have \(\Psi_\delta(x)=\sum_{i=1}^d (x_i^2+\delta^2)^{q/2}\le\sum_{i=1}^d |x_i|^q+d\delta^q\le D_X^q+d\delta^q ,\) which derives
\(    \Lambda_{ii}(x)
    \le
    \frac{1}{q-1}
    (D_X^q+d\delta^q)^{(2-q)/q}\delta^{q-2}, \forall i=1,\ldots,d\) and \(x\in X.
\) 
Since $\|y\|_2\le \|y\|_q$ for $q\in(1,2]$, we obtain \[y^\top \Lambda(x)y\le\frac{1}{q-1}(D_X^q+d\delta^q)^{(2-q)/q}\delta^{q-2}\|y\|_q^2 .\]
Combining the bounds for $y^\top M(x)y$ and $y^\top \Lambda(x)y$, we have
\(
    y^\top \nabla^2 v_\delta(x)y
    \le
    L_{v}\|y\|_q^2,
    \) \(\forall x\in X,  y\in\mathbb R^d,
\)
where $L_v$ is given by \eqref{eq:Lv}.
We next show that this quadratic-form bound implies the corresponding
operator-norm bound. Since $\nabla^2 v_\delta(x)$ is symmetric positive
semidefinite, the Cauchy--Schwarz inequality w.r.t. the
semidefinite bilinear form induced by $\nabla^2 v_\delta(x)$ gives, for
any $y,z\in\mathbb{R}^d$,
\[
    |z^\top \nabla^2 v_\delta(x) y|
    \le
    \sqrt{z^\top \nabla^2 v_\delta(x) z}\,
    \sqrt{y^\top \nabla^2 v_\delta(x) y} 
    \le
    L_v \|z\|_q \|y\|_q .
\]
Then by the duality between the $\ell_q$ and $\ell_p$ norms,
\(\|\nabla^2 v_\delta(x)y\|_p
    = \sup_{\|z\|_q\le 1}
       |z^\top \nabla^2 v_\delta(x)y| 
    \le L_v\|y\|_q .\)
Hence,
\[
    \|\nabla^2 v_\delta(x)\|_p
    =\sup_{\|y\|_q\le1}
      \|\nabla^2 v_\delta(x)y\|_p
    \le L_v,
    \qquad \forall x\in X.
\]
Since $X$ is convex, for any $x,y\in X$, the segment $x+t(y-x)$ remains in $X$. Hence,
\[\|\nabla v_\delta(y)-\nabla v_\delta(x)\|_p\le\int_0^1\|\nabla^2v_\delta(x+t(y-x))(y-x)\|_p\,dt \le L_{v}\|y-x\|_q .\]
Therefore, $v_\delta$ is $L_{v}$-smooth on $X$ w.r.t. the $\|\cdot\|_q$-norm, which completes the proof.
\end{proof}

\section{Proofs of Theorem \ref{thm:restart} and Corollary \ref{cor:2-restart}}
\label{app:restart}

\begin{proof}[Proof of Theorem~\ref{thm:restart}]
For \(s\)-th stage, let \(\mathcal F_{s-1}\) denote the \(\sigma\)-field generated by all randomness before stage \(s\), and define \(\Delta_{0,s}=\mathcal L(x_s^0,0;\mu_s)-\Phi_{\rm inf}+2m\rho \bar{M}_c^2\). By assumption, the regularity condition in Assumption~\ref{ass:ensc} holds at every iterate generated during stage \(s\), with the same positive constant \(\beta\) for all stages. Hence, conditioned on \(\mathcal F_{s-1}\), the stationarity and feasibility bounds established for Algorithm~\ref{algo} can be applied to each stage with the corresponding stagewise parameters.

We first verify that \(\mu_s\) and \(\epsilon_s\) are compatible with the parameter requirements in \eqref{settings}. By the definitions of \(\gamma\) and \(\epsilon_s\), \(40(M_f^2+M_h^2+m^2\rho^2\bar{M}_c^2M_c^2)/(\beta^2\epsilon_s^2)=\gamma^2/\epsilon_s^2=\mu_s^2\). Moreover, since \(\mu_s\ge\mu_{\rm base}\), we have \(\mu_s^2\ge 10(8L_v^2+3)/(27\beta^2)\) and \(\mu_s^2\ge8(2L_v^2+1)/(3\beta^2)\). Hence, the penalty requirement in \eqref{settings} is satisfied at accuracy \(\epsilon_s\). Also, \(\mu_s\ge\mu_{\rm base}\ge\gamma\), so \(0<\epsilon_s\le1\). All remaining batch-size, smoothing, momentum, and stepsize parameters are chosen according to \eqref{settings}, with \(\epsilon\) and \(\mu\) replaced by \(\epsilon_s\) and \(\mu_s\), respectively; when \(p>2\ln d\), the replacements specified in Theorem~\ref{thm:pinf} are used. 
For stage \(s=1\), since \(x_1^0=x^0\), Assumption~\ref{ass:lip} gives \(\|c(x_1^0)\|_2^2\le m\bar{M}_c^2=\bar C_1\), while \(\Delta_{0,1}\le B_X+\mu_1m\bar{M}_c^2/2=\bar\Delta_1\). Hence, \(\mathbb E[\Delta_{0,1}]\le\bar\Delta_1\) and \(\mathbb E\|c(x_1^0)\|_2^2\le\bar C_1\). We now prove inductively that, for every stage \(s=1,\ldots,T,\)
\begin{equation}\label{eq:restart-initial-bounds}
\mathbb E[\Delta_{0,s}]\le\bar\Delta_s,
\qquad
\mathbb E\|c(x_s^0)\|_2^2\le\bar C_s,
\end{equation}
and that the output of stage \(s\) satisfies
\begin{align}
\mathbb E[\|\mathcal G_X(\tilde x^s,\nabla f(\tilde x^s)+\nabla c(\tilde x^s)^\top\tilde\lambda^s,\eta_s)\|_q^2]
&\le\epsilon_s^2,
\label{eq:restart-stat}\\
\mathbb E\|c(\tilde x^s)\|_p^2
&\le\epsilon_s^2,
\label{eq:restart-feas}
\end{align}
where \(\tilde\lambda^s\) is the multiplier associated with the random output of stage \(s\).

Suppose that \eqref{eq:restart-initial-bounds} holds at \(s\)-th stage with $s\ge1$. Conditioned on \(\mathcal F_{s-1}\), the starting point \(x_s^0\) is fixed, so the stationarity bound \eqref{bound} applies conditionally. Taking total expectation, the first three terms on the right-hand side of \eqref{bound} are bounded by \(\epsilon_s^2/4\) each by the same parameter substitutions as in the proofs of Theorems~\ref{thm:p} and \ref{thm:pinf}, while the last term satisfies \(27\mathbb E[\Delta_{0,s}]/(\eta_sK_s)\le27\bar\Delta_s/(\eta_sK_s)\le\epsilon_s^2/4\). Therefore, \eqref{eq:restart-stat} holds. When \(p>2\ln d\), the stochastic estimation terms are controlled in the \(\ell_{\bar p}\)-norm with \(\bar p=2\ln d\), exactly as in the proof of Theorem~\ref{thm:pinf}.
For primal feasibility, conditioned on \(\mathcal F_{s-1}\), the bound \eqref{bound-fea} applies. Using \(\|z\|_p\le\|z\|_2\) for \(p\ge2\), its first term satisfies \(\mathbb E\|c(x_s^0)\|_p^2/K_s\le\bar C_s/K_s\le\epsilon_s^2/5\). For the term involving \(\Delta_{0,s}\), it follows from  \(\mu_s^2\ge10(8L_v^2+3)/(27\beta^2)\) that 
\[
\frac{8(8L_v^2+3)\mathbb E[\Delta_{0,s}]}{\beta^2\mu_s^2\eta_sK_s}
\le
\frac{2(8L_v^2+3)}{27\beta^2\mu_s^2}\epsilon_s^2
\le
\frac{\epsilon_s^2}{5}.
\]
 The remaining three terms in \eqref{bound-fea} are bounded by \(\epsilon_s^2/5\) each by the same parameter substitutions as in the proofs of Theorems~\ref{thm:p} and \ref{thm:pinf}. Consequently, \eqref{eq:restart-feas} holds.
It remains to verify that the output of stage \(s\) provides the required initial bounds for stage \(s+1\). Since \(x_{s+1}^0=\tilde x^s\)  and \(p\ge2\), \(\|c(x)\|_2^2\le m^{1-2/p}\|c(x)\|_p^2\). Therefore, by \eqref{eq:restart-feas}, \(\mathbb E\|c(x_{s+1}^0)\|_2^2\le m^{1-2/p}\epsilon_s^2=m^{1-2/p}\gamma^2/\mu_s^2\). Since \(\mu_{s+1}\le2\mu_s^2\), we obtain 
\(%begin{equation}\label{eq:restart-near-feas}
\mathbb E\|c(x_{s+1}^0)\|_2^2
\le\frac{2m^{1-2/p}\gamma^2}{\mu_{s+1}}
=\bar C_{s+1}.
\) %end{equation}
By the definition of \(B_X\) we further obtain 
\(\mathbb E[\Delta_{0,s+1}]
\le B_X+\frac{\mu_{s+1}}{2}\mathbb E\|c(x_{s+1}^0)\|_2^2
\le B_X+m^{1-2/p}\gamma^2
=\bar\Delta_+.\)
Since \(s+1\ge2\), \(\bar\Delta_{s+1}=\bar\Delta_+\), so \(\mathbb E[\Delta_{0,s+1}]\le\bar\Delta_{s+1}\). Thus, \eqref{eq:restart-initial-bounds} holds at stage \(s+1\), and the induction is complete.

By construction, \(\mu_T\ge \mu_\star\ge\gamma/\epsilon\), hence \(\epsilon_T=\gamma/\mu_T\le\epsilon\). Combining this with \eqref{eq:restart-stat}--\eqref{eq:restart-feas} shows that \(\tilde x^T\) is an \(\epsilon\)-KKT point. 
% To bound the number of  stages, if \(\mu_\star=\mu_1=\mu_{\rm base}\), then \(T=1\). Otherwise, before the  \(\mu_\star\) is reached, \(\mu_{s+1}=2\mu_s^2\). Setting \(y_s=2\mu_s\) gives \(y_{s+1}=y_s^2\), it derives \(\mu_s=\frac12(2\mu_1)^{2^{s-1}}\). Consequently,
% \[
% T\le1+\lceil\log_2(\frac{\log(2\mu_\star)}{\log(2\mu_1)})\rceil
% =\mathcal O(1+\log\log\frac{e}{\epsilon}).\]
The remainder is to bound the total number of stochastic function evaluations, treating the first stage separately from the later stages. By the definition of \(\mu_{\rm base}\), the expression \eqref{eq:Lv} and the convention on problem constants used in Theorem~\ref{thm:rad_qnorm}, we obtain that \(\mu_s\ge\mu_{\rm base}\ge\sqrt{8(2L_v^2+1)/(3\beta^2)}\) for any $s\ge1$, so \(L_v=\mathcal O(\mu_s)=\mathcal O(\epsilon_s^{-1})\), where \(\mu_s=\gamma/\epsilon_s\). By the same argument as in the proof of Theorem~\ref{thm:rad_qnorm}, \(\eta_s^{-1}=\mathcal O(\epsilon_s^{-1}+L_v)=\mathcal O(\epsilon_s^{-1})\). Therefore, 
for \(s=1\), we have \(\bar C_1=m\bar{M}_c^2\), \(\epsilon_1^{-1}=\mu_{\rm base}/\gamma=\mathcal O(1+L_v)\), and \(\eta_1^{-1}=\mathcal O(\epsilon_1^{-1})=\mathcal O(1+L_v)\). Hence
\[
K_1=\mathcal O(\frac{\bar\Delta_1}{\eta_1\epsilon_1^2}
+\frac{\bar C_1}{\epsilon_1^2}+\eta_1^{-2})
=\mathcal O((1+L_v)^4).
\]
For \(s\ge2\), \(\bar C_s/\epsilon_s^2=2m^{1-2/p}\mu_s=\mathcal O(\epsilon_s^{-1})\), and then 
\[
K_s=\mathcal O(\frac{\bar\Delta_s}{\eta_s\epsilon_s^2}
+\frac{\bar C_s}{\epsilon_s^2}+\eta_s^{-2})
=\mathcal O(\epsilon_s^{-3}),
\qquad s\ge2.
\]
To count the total number of stochastic function value evaluations, we consider \(p\in[2,2\ln d]\) and \(p>2\ln d\), respectively. For $p\in[2,2\ln d]$, under Rademacher smoothing, \(n_{0,s}=\mathcal O(pd^{2/p}\epsilon_s^{-2})\) and \(n_s=\mathcal O(pd^{2/p})\). Hence, for every \(s\ge2\),  \(2n_{0,s}+4n_sK_s=\mathcal O(pd^{2/p}\epsilon_s^{-3})\). For \(s=1\), using \(L_v=\mathcal O(pd^{1-2/p})\) and therefore \(1+L_v=\mathcal O(pd^{1-2/p})\) for \(p\ge2\) and \(d\ge3\),
\[
2n_{0,1}+4n_1K_1
=\mathcal O(pd^{2/p}(1+L_v)^4)
=\mathcal O(p^5d^{4-6/p}).\]  
For \(p>2\ln d\), Theorem~\ref{thm:pinf} gives \(n_{0,s}=\mathcal O(\ln d\cdot\epsilon_s^{-2})\) and \(n_s=\mathcal O(\ln d)\). Thus, for \(s\ge2\), \(2n_{0,s}+4n_sK_s=\mathcal O(\ln d\cdot\epsilon_s^{-3})\), while the first stage satisfies
\[
2n_{0,1}+4n_1K_1
=\mathcal O((1+L_v)^4\ln d)
=\mathcal O(p^4d^{4-8/p}\ln d).\]
Therefore, the total number of stochastic function value evaluations is 
\be\label{sum-cp}
2n_{0,1} + 4n_1K_1 + \sum_{s= 2}^T (2n_{0,s} + 4n_s K_s) = \begin{cases}
\mathcal O(p^5d^{4-6/p}) + \mathcal O(pd^{2/p} \sum_{s=2}^T\epsilon_s^{-3} ),& \quad p\in[2,\ln d],\\
\mathcal O(p^4d^{4-8/p}\ln d) + \mathcal O(\ln d\cdot \sum_{s=2}^T \epsilon_s^{-3}),& \quad p>2\ln d.
\end{cases}
\ee
It remains to sum the oracle costs of the stages \(s\ge2\). If \(T=1\), there is nothing to sum. Suppose \(T\ge2\). Then \(\mu_\star>\mu_{\rm base}\), so \(\mu_T=\mu_\star=\gamma/\epsilon\) and \(\epsilon_T=\epsilon\). For every stage before the final capped stage, i.e., \(s+1<T\), since \(\mu_s\ge1\), we have \(\mu_{s+1}=2\mu_s^2\ge2\mu_s\) and hence \(\epsilon_{s+1}^{-3}\ge8\epsilon_s^{-3}\). 
At the final capped stage, \(\mu_T\ge\mu_{T-1}\), so \(\epsilon_{T-1}^{-3}\le\epsilon_T^{-3}\).   Therefore, we obtain
\[
\sum_{s=2}^{T}\epsilon_s^{-3} = \sum_{s=2}^{T-1}\epsilon_s^{-3} + {\epsilon_T}^{-3}
\le(1+\frac18+\frac1{8^2}+\cdots+\frac{1}{8^{T-3}})\epsilon_{T-1}^{-3}
 + \epsilon_T^{-3} \le\frac87\epsilon_{T-1}^{-3}  + \epsilon_T^{-3}=\mathcal O(\epsilon^{-3}).
\]
Then combining with \eqref{sum-cp} yields the final oracle complexity bounds \eqref{mul-cp}.
% the first-stage cost with the sum over \(s\ge2\), the total oracle complexity is \(\mathcal O(pd^{2/p}\epsilon^{-3}+p^5d^{4-6/p})\) for \(p\in[2,2\ln d]\), and \(\mathcal O(\ln d\cdot\epsilon^{-3}+p^4d^{4-8/p}\ln d)\) for \(p>2\ln d\). This proves the claimed oracle complexity bounds.
\end{proof}

\begin{proof}[Proof of Corollary~\ref{cor:2-restart}]
The result follows by specializing the multi-stage analysis of
Theorem~\ref{thm:restart} to the Euclidean case \(p=2\). 
In this case, \(v(x)=\frac12\|x\|_2^2,\) \(L_v=1,\) \(S_2=d .\)
Thus, at the first stage,  \(\epsilon_1^{-1}=\frac{\mu_{\rm base}}{\gamma}=\mathcal O(1+L_v)=\mathcal O(1),\) and \(\eta_1^{-1}=\mathcal O(\epsilon_1^{-1}+L_v)=\mathcal O(1).\)  It derives \(K_1=\mathcal O(1)\).
Moreover, the total number of stochastic function evaluations of the first stage is
\( 
2n_{0,1}+4n_1K_1 =\mathcal O(d).
\) 
For every stage \(s\ge2\), Corollary \ref{thm:2-norm}
with parameter settings \(\epsilon\) replaced by \(\epsilon_s\) and \(\mu\) replaced by
\(\mu_s\) gives
\(
    2n_{0,s}+4n_sK_s
    =
    \mathcal O(d\epsilon_s^{-3}).
\) 
Therefore, it follows from \(0<\epsilon\le1\) that the total oracle complexity is
\( \mathcal O(d)+\mathcal O(d\epsilon^{-3})=\mathcal O(d\epsilon^{-3})\).
\end{proof}
\end{document}